\documentclass{amsproc}

\usepackage{amsmath,amssymb}
\newtheorem{theorem}{Theorem}[section]
\newtheorem{lemma}[theorem]{Lemma}
\newtheorem{prop}[theorem]{Proposition}

\theoremstyle{definition}
\newtheorem{definition}[theorem]{Definition}

\newtheorem{coro}[theorem]{Corollary}
\theoremstyle{remark}
\newtheorem{remark}[theorem]{Remark}

\numberwithin{equation}{section}

\def\lg{\langle}
\def\rg{\rangle}

\def\lra{\longrightarrow}
\def\al{\alpha}
\def\be{\beta}

\def\si{\sigma}

\def\vn{\varepsilon}

\def\ep{\epsilon}
\def\ot{\otimes}

\def\om{\omega}
\def\lra{\longrightarrow}

\begin{document}

\title[Representations of two-parameter quantum groups]
{Representations of Two-parameter Quantum Orthogonal and
Symplectic Groups}

\author[Bergeron]{Nantel Bergeron}
\address{Department of Mathematics and Statistics, York
University, North York, Toronto, ON Canada, M3J 1P3}
\email{bergeron@mathstat.yorku.ca}
\thanks{N.B. supported by the NSERC and the CRC of
Canada.}

\author[Gao]{Yun Gao}
\address{Department of Mathematics and Statistics, York
University, North York, Toronto, ON Canada, M3J 1P3}
\thanks{Y.G. supported by the NSERC of Canada and the Chinese Academy of Science.}
\email{ygao@yorku.ca}

\author[Hu]{Naihong Hu$^\star$}
\address{Department of Mathematics, East China Normal University,
Shanghai 200062, PR China} \email{nhhu@euler.math.ecnu.edu.cn}
\thanks{$^\star$N.H., Corresponding Author,
supported in part by the NNSF (Grant 10431040), the TRAPOYT and the
FUDP from the MOE of China, the SRSTP from the STCSM, the Shanghai
Priority Academic Discipline from the SMEC}

\date{Sept. 1, 2004 and, in revised form, May 18, 2005.}


\begin{abstract}
We investigate the finite-dimensional representation theory of
two-parameter quantum orthogonal and symplectic groups that we found
in [BGH] under the assumption that $rs^{-1}$ is not a root of unity
and extend some results [BW1, BW2] obtained for type $A$ to types
$B$, $C$ and $D$. We construct the corresponding $R$-matrices and
the quantum Casimir operators, by which we prove that the complete
reducibility Theorem also holds for the categories of
finite-dimensional weight modules for types $B$, $C$, $D$.
\end{abstract}

\maketitle



\section{Preliminaries: Two-parameter Quantum Groups for Classical Types} 

Let ${\mathbb K}\supset{\mathbb Q}(r,s)$ denote an algebraically
closed field, where the two-parameters $r,\,s$ are nonzero complex
numbers satisfying $r^2\ne s^2$.

In this section, we recall the definitions of the two-parameter
quantum groups $U_{r,s}({\mathfrak g})$ for ${\mathfrak
g}={\mathfrak{sl}}_{n+1}$ from [BW1], and for ${\mathfrak
g}={\mathfrak{so}}_{2n+1}$, ${\mathfrak{sp}}_{2n}$ and ${\mathfrak
{so}}_{2n}$ from [BGH]. Let $\Psi$ be a finite root system of a
simple Lie algebra ${\mathfrak g}$ of rank $n$ with $\Pi$ a base of
simple roots. Regard $\Psi$ as a subset of a Euclidean space
$E={\mathbb R}^n$ with an inner product $(\,,)$. Let $\ep_1, \cdots,
\ep_n$ denote an orthonormal basis of $E$. We need the following
data on (prime) root systems.

Type $A$:
\begin{equation*}
\begin{split}
\Pi&=\{\al_i=\ep_i-\ep_{i+1}\mid 1\le i\le n\},\\
\Psi&=\{\pm(\ep_i-\ep_j)\mid 1\le i<j\le n+1\}.
\end{split}
\end{equation*}

Type $B$:
\begin{equation*}
\begin{split}
\Pi&=\{\al_i=\ep_i-\ep_{i+1}\mid 1\le i<n\}\cup\{\al_n=\ep_n\},\\
\Psi&=\{\pm\ep_i\pm\ep_j\mid 1\le i\ne j\le n\}\cup\{\pm\ep_i\mid
1\le i\le n\}.
\end{split}
\end{equation*}

Type $C$:
\begin{equation*}
\begin{split}
\Pi&=\{\al_i=\ep_i-\ep_{i+1}\mid 1\le i<n\}\cup\{\al_n=2\ep_n\},\\
\Psi&=\{\pm\ep_i\pm\ep_j\mid 1\le i\ne j\le n\}\cup\{2\ep_i\mid
1\le i\le n\}.
\end{split}
\end{equation*}

Type $D$:
\begin{equation*}
\begin{split}
\Pi&=\{\al_i=\ep_i-\ep_{i+1}\mid 1\le
i<n\}\cup\{\al_n=\ep_{n-1}+\ep_n\},\\
\Psi&=\{\pm\ep_i\pm\ep_j\mid 1\le i\ne j\le n\}.
\end{split}
\end{equation*}

In the cases of type $A$, $C$ and $D$, we set
$r_i=r^{\frac{(\al_i,\al_i)}2}$, $s_i=s^{\frac{(\al_i,\al_i)}2}$;
while for type $B$, we set $r_i=r^{(\al_i, \al_i)}$, $s_i=s^{(\al_i,
\al_i)}$.

Assigned to $\Pi$, there are two sets of mutually-commutative
symbols $W=\{\om_i^{\pm1}\mid 1\le i\le n\}$ and
$W'=\{{\om_i'}^{\pm1}\mid 1\le i\le n\}$. Define a pairing $\lg\,
, \rg:\, W'\times W\lra {\mathbb K}$ as follows:

\begin{equation*}
\lg \om_i', \om_j\rg= r^{(\ep_j,\al_i)}s^{(\ep_{j+1},\al_i)},
\quad\qquad  i\le n{+}1, \ j\le n, \ \qquad\;\text{\it for}\quad
\mathfrak{sl}_{n+1},\leqno(\text{$1_A$})
\end{equation*}
\begin{equation*}
\lg \om_i', \om_j\rg=\begin{cases}
r^{2(\ep_j,\al_i)}s^{2(\ep_{j+1},\al_i)}, & \quad\qquad\ i\le n, \
j<n,\cr r^{2(\ep_n, \al_i)}, & \quad\qquad\ i<n, \ j=n,\qquad\
\text{\it for}\quad \mathfrak{so}_{2n+1},\cr r^{(\ep_n,
\al_n)}s^{-(\ep_n, \al_n)}, & \qquad\quad\   i=j=n.
\end{cases}\tag{$1_B$}
\end{equation*}
\begin{equation*}
\lg \om_i', \om_j\rg=\begin{cases}
r^{(\ep_j,\al_i)}s^{(\ep_{j+1},\al_i)}, & \quad\qquad  i\le n, \
j<n,\cr r^{2(\ep_n, \al_i)}, & \quad\qquad i<n, \ j=n, \
\qquad\;\,\text{\it for}\quad \mathfrak{sp}_{2n},\cr r^{(\ep_n,
\al_n)}s^{-(\ep_n, \al_n)}, & \quad\qquad i=j=n.
\end{cases}\tag{$1_C$}
\end{equation*}
\begin{equation*}
\lg \om_i', \om_j\rg=\begin{cases}
r^{(\ep_j,\al_i)}s^{(\ep_{j+1},\al_i)}, & \ \;  i\le n, \ j<n,\cr
r^{(\ep_{n-1}, \al_i)}s^{-(\ep_n,\al_i)}, & \ \; i\ne n-1, \
j=n,\quad\text{\it for}\quad \mathfrak{so}_{2n},\cr r^{(\ep_n,
\al_{n-1})}s^{-(\ep_{n-1}, \al_{n-1})}, & \ \;i=n-1,\ j=n.
\end{cases}\tag{$1_D$}
\end{equation*}
\begin{equation*}\lg {\om_i'}^{\pm1},\om_j^{-1}\rg=\lg
{\om_i'}^{\pm1},\om_j\rg^{-1}=\lg \om_i',\om_j\rg^{\mp1},\qquad
\text{\it for any } \ \mathfrak g.  \tag{2}
\end{equation*}

\smallskip

\begin{lemma} \ For the prime root systems of the Lie algebras $\mathfrak g=\mathfrak
{sl}_{n}$, $\mathfrak{so}_{2n+1}$, $\mathfrak{so}_{2n}$, and
$\mathfrak{sp}_{2n}$, there hold the identities:
\begin{gather*}
(\ep_{j+1},\al_i)=-(\ep_i,\al_j),\, \qquad (i, j<n),
\;\qquad\qquad\quad
\text{\it for any } \ \mathfrak g,\\
(\ep_{j+1},\al_n)=\begin{cases}-(\ep_n,\al_j), \quad &\quad (j<n),
\qquad\quad\quad \qquad \text{\it
for } \ \mathfrak g=\mathfrak{so}_{2n+1}, \\
-2(\ep_n,\al_j), \quad &\quad (j<n), \quad \qquad\qquad\quad
\text{\it for } \ \mathfrak g=\mathfrak{sp}_{2n},
\end{cases}\\
(\ep_j,\al_n)=\begin{cases}-(\ep_n,\al_{j-1}),  &\quad (j\le n,\, j\ne n-1),  \\
(\ep_{n-1},\al_{n-1}), &\quad (j=n-1)
\end{cases}
\,\quad\text{\it for } \ \mathfrak g=\mathfrak{so}_{2n}.
\end{gather*}\hfill\qed
\end{lemma}

Observe that Lemma 1.1 ensures the compatibility of the defining
relations of the two-parameter quantum groups defined below.

Let $U_{r,s}(\mathfrak g)$ be the unital associative algebra over
$\mathbb K$ generated by symbols $e_i, f_i, \omega_i^{\pm 1}$,
${\om_i'}^{\pm 1}$ $(1\le i\le n)$, subject to the following
relations $(X1)$---$(X4)$:

\smallskip
$(X1)$ \  $\om_i^{\pm
1}{\om_j'}^{\pm1}={\om_j'}^{\pm1}{\om_i}^{\pm1}$,\qquad
$\om_i^{\pm 1}{\om_i}^{\mp1}=1={\om_i'}^{\pm1}{\om_i'}^{\mp1}$.

\smallskip
$(X2)$ \ For $1\le i, \,j \le n$, we have
\begin{alignat*}{2}
    \om_j\,e_i\,\om_j^{-1}&=\lg \om_i',\om_j\rg\,e_i,  &\hspace{20pt}
    \om_j\,f_i\,\om_j^{-1}&=\lg \om_i',\om_j\rg^{-1}\,f_i,\\
    \om_j'\,e_i\,{\om_j'}^{-1}&=\lg \om_j',\om_i\rg^{-1}\,e_i,
   &\hspace{20pt} \om_j'\,f_i\,{\om_j'}^{-1}&=\lg
   \om_j',\om_i\rg\,f_i.\label{eq:x}
\end{alignat*}

$(X3)$ \ For $1\le i,\, j\le n$, we have
$$
[\,e_i, f_j\,]=\delta_{ij}\frac{\om_i-\om_i'}{r_i-s_i}.
$$

\smallskip
$(X4)$ \ For any $i\ne j$, we have the $(r,s)$-Serre relations:
\begin{gather*}
\bigl(\text{ad}_l\,e_i\bigr)^{1-a_{ij}}\,(e_j)=0,\\
\bigl(\text{ad}_r\,f_i\bigr)^{1-a_{ij}}\,(f_j)=0,
\end{gather*}
where the definitions of the left-adjoint action
$\text{ad}_l\,e_i$ and the right-adjoint action $\text{ad}_r\,f_i$
are given in the following sense:
$$
\text{ad}_{ l}\,a\,(b)=\sum_{(a)}a_{(1)}\,b\,S(a_{(2)}), \quad
\text{ad}_{ r}\,a\,(b)=\sum_{(a)}S(a_{(1)})\,b\,a_{(2)}, \quad
\forall\; a, b\in U_{r,s}(\mathfrak g),
$$
where $\Delta(a)=\sum_{(a)}a_{(1)}\ot a_{(2)}$ is given by
Proposition 1.2 below.

\smallskip
The following fact is straightforward.
\begin{prop}  \ The algebra $U_{r, s}(\mathfrak g)$
$($\,$\mathfrak g=\mathfrak{sl}_{n+1}$,
$\mathfrak{so}_{2n+1},\,\mathfrak{sp}_{2n}$, or
$\mathfrak{so}_{2n}$\,$)$ is a Hopf algebra under the
comultiplication, the counit and the antipode defined below:
\begin{gather*}
\Delta(\om_i^{\pm1})=\om_i^{\pm1}\ot\om_i^{\pm1},
\qquad
\Delta({\om_i'}^{\pm1})={\om_i'}^{\pm1}\ot{\om_i'}^{\pm1},\\
\Delta(e_i)=e_i\ot 1+\om_i\ot e_i, \qquad \Delta(f_i)=1\ot
f_i+f_i\ot \om_i',\\
\vn(\om_i^{\pm})=\vn({\om_i'}^{\pm1})=1, \qquad
\vn(e_i)=\vn(f_i)=0,\\
S(\om_i^{\pm1})=\om_i^{\mp1}, \qquad
S({\om_i'}^{\pm1})={\om_i'}^{\mp1},\\
S(e_i)=-\om_i^{-1}e_i,\qquad S(f_i)=-f_i\,{\om_i'}^{-1}.
\end{gather*}\hfill\qed
\end{prop}

\begin{remark}
\ When $r=s^{-1}=q$, Hopf algebra $U_{r, s}(\mathfrak g)$ modulo
the Hopf ideal generated by the elements $\om_i'-\om_i^{-1}$
$(1\le i\le n)$, is just the quantum groups $U_q(\mathfrak g)$ of
Drinfel'd-Jimbo type.
\end{remark}

\begin{definition} \ {\it A skew-dual pairing of two Hopf
algebras ${\mathcal A}$ and ${\mathcal U}$ is a bilinear form
$\langle\,,\rangle:\; {\mathcal U}\times {\mathcal A}\lra \mathbb
K$ such that
\begin{gather*} \langle f, 1_{\mathcal A}\rangle=\vn_{\mathcal U}(f),\qquad
\langle 1_{\mathcal U}, a\rangle=\vn_{\mathcal A}(a),\\
\langle f, a_1a_2\rangle=\langle \Delta^{\text{op}}_{\mathcal
U}(f), a_1\ot a_2\rangle, \qquad \langle f_1f_2, a \rangle=\langle
f_1\ot f_2, \Delta_{\mathcal A}(a)\rangle,
\end{gather*}
for all $f,\, f_1,\, f_2\in\mathcal U$, and
$a,\,a_1,\,a_2\in\mathcal A$, where $\vn_{\mathcal U}$ and
$\vn_{\mathcal A}$ denote the counits of $\mathcal U$ and
$\mathcal A$, respectively, and $\Delta_{\mathcal U}$ and
$\Delta_{\mathcal A}$ are their respective comultiplications.}
\end{definition}

Let $\mathcal B=B(\mathfrak g)$ (resp. $\mathcal B'=B'(\mathfrak
g)$\,) denote the Hopf subalgebra of $U=U_{r,s}(\mathfrak g)$
generated by $e_j$, $\om_j^{\pm1}$ (resp. $f_j$,
${\om_j'}^{\pm1}$\,) with $1\le j\le n$ for $\mathfrak
g=\mathfrak{sl}_{n+1}$, and with $1\le j\le n$ for $\mathfrak
g=\mathfrak {so}_{2n+1}$, $\mathfrak {so}_{2n}$, and
$\mathfrak{sp}_{2n}$, respectively. The following result was
obtained for the type $A$ case by [BW1], and for the types $B$, $C$
and $D$ cases by [BGH].

\begin{prop} \ There exists a unique skew-dual pairing
$\langle\,,\rangle:\, \mathcal B'\times \mathcal B\lra\mathbb K$
of the Hopf subalgebras $\mathcal B$ and $\mathcal B'$ in
$U_{r,s}(\mathfrak g)$, for $\mathfrak g=\mathfrak {sl}_{n+1}$,
$\mathfrak {so}_{2n+1}$, $\mathfrak {so}_{2n}$, or
$\mathfrak{sp}_{2n}$ such that $\lg f_i,
e_j\rg=\frac{\delta_{ij}}{s_i-r_i}$, and the conditions $(1_X)$
$($where $X=A,\,B,\,C$, or $D$$)$ and $(2)$ are satisfied, and all
other pairs of generators are $0$. Moreover, we have $\lg S(a),
S(b)\rg=\lg a, b\rg$ for $a\in\mathcal B',\,b\in\mathcal
B$.\hfill\qed
\end{prop}

\begin{definition} \ {\it For any two skew-paired Hopf algebras
$\mathcal A$ and $\mathcal U$ by a skew-dual pairing $\lg\,,\rg$,
one may form the Drinfel'd double $\mathcal D(\mathcal A,\mathcal
U)$ as in $[\rm KS, 8.2]$, which is a Hopf algebra whose
underlying coalgebra is $\mathcal A\ot\mathcal U$ with the tensor
product coalgebra structure, and whose algebra structure is
defined by
$$
(a\ot f)(a'\ot f')=\sum \lg S_{\mathcal U}(f_{(1)}),
a'_{(1)}\rg\lg f_{(3)},a'_{(3)}\rg \,aa'_{(2)}\ot
f_{(2)}f',\leqno(3)$$ for $a, a'\in \mathcal A$ and $f,
f'\in\mathcal U$. The antipode $S$ is given by}
$$
S(a\ot f)=(1\ot S_{\mathcal U}(f))(S_{\mathcal A}(a)\ot 1).
$$
\end{definition}

Clearly, both mappings $\mathcal A\ni a\mapsto a\ot 1\in\mathcal
D(\mathcal A,\mathcal U)$ and $\mathcal U\ni f\mapsto 1\ot
f\in\mathcal D(\mathcal A, \mathcal U)$ are injective Hopf algebra
homomorphisms. Let us denote the image $a\ot 1$ (resp. $1\ot f$)
of $a$ (resp. $f$) in $\mathcal D(\mathcal A,\mathcal U)$ by $\hat
a$ (resp. $\hat f$). By (3), we have the following cross
commutation relations between elements $\hat a$ (for $a\in\mathcal
A$) and $\hat f$ (for $f\in\mathcal U$) in the algebra $\mathcal
D(\mathcal A,\mathcal U)$:
\begin{gather*} \hat f\,\hat a=\sum\, \lg
S_{\mathcal U}(f_{(1)}), a_{(1)}\rg\,\lg
f_{(3)},a_{(3)}\rg\;\hat a_{(2)}\hat f_{(2)},\tag{4}\\
\sum\lg f_{(1)}, a_{(1)}\rg\,\hat f_{(2)}\,\hat a_{(2)}= \sum \hat
a_{(1)}\,\hat f_{(1)}\,\lg f_{(2)},a_{(2)}\rg.\tag{5}
\end{gather*}
In fact, as an algebra the double $\mathcal D(\mathcal A,\mathcal
U)$ is the universal algebra generated by the algebras $\mathcal
A$ and $\mathcal U$ with cross relations (4) or, equivalently,
(5).

\begin{theorem}[{[BW1, BGH]}] \ The two-parameter quantum group
$U=U_{r,s}(\mathfrak g)$  is isomorphic to the Drinfel'd quantum
double $\mathcal D(\mathcal B, \mathcal B')$, for $\mathfrak
g=\mathfrak {sl}_{n+1}$, $\mathfrak {so}_{2n+1}$, $\mathfrak
{so}_{2n}$, or $\mathfrak{sp}_{2n}$.\hfill\qed
\end{theorem}

Let us denote $U_{r,s}(\mathfrak n)$ $($resp. $U_{r,s}(\mathfrak
n^-)$\,$)$ the subalgebra of $\mathcal B$ $($resp. $\mathcal
B'$$)$ generated by $e_i$ $($resp. $f_i$$)$ for all $i\le n$. Let
\begin{gather*}
U^0=\mathbb
K\,[\om_1^{\pm1},\cdots,\om_n^{\pm1},{\om_1'}^{\pm1},\cdots,{\om_n'}^{\pm1}],\\
U_0=\mathbb K\,[\om_1^{\pm1},\cdots,\om_n^{\pm1}],\qquad
U_0'=\mathbb K\,[{\om_1'}^{\pm1},\cdots,{\om_n'}^{\pm1}]
\end{gather*}
denote the respective Laurent polynomial subalgebras of
$U_{r,s}(\mathfrak g)$, $\mathcal B$,
 and $\mathcal B'$. Clearly, $U^0=U_0U_0'=U_0'U_0$. Thus, by definition, we have
$\mathcal B=U_{r,s}(\mathfrak n)\rtimes U_0$, and $\mathcal
B'=U_0'\ltimes U_{r,s}(\mathfrak n^-)$, such that the double
$\mathcal D(\mathcal B,\mathcal B')\cong U_{r,s}(\mathfrak n)\ot
U^0\ot U_{r,s}(\mathfrak n^-)$, as vector spaces.

Let $\lg\,|\,\rg_0:\,\mathcal B\times \mathcal B'\lra \mathbb K$
denote the skew-dual pairing given by $\lg b\,|\,b'\rg_0=\lg S(b'),
b\rg$. Then, via a variation of its Drinfel'd double structure, we
obtain the standard triangular decomposition of $U_{r,s}(\mathfrak
g)$ in [BGH, Corollary 2.6] as follows.

\begin{coro} \ $U_{r,s}(\mathfrak g)\cong U_{r,s}(\mathfrak
n^-)\ot U^0\ot U_{r,s}(\mathfrak n)$, {\it as vector spaces. In
particular, it induces $U_q(\mathfrak g)\cong U_q(\mathfrak
n^-)\otimes U_0\otimes U_q(\mathfrak n)$, as vector
spaces.\hfill\qed}
\end{coro}

Let $Q=\mathbb Z\Psi$ denote the root lattice and set
$Q^+=\sum_{i=1}^n\mathbb Z_{\ge0}\al_i$. Then for any
$\zeta=\sum_{i=1}^n\zeta_i\al_i\in Q$, we denote
$$
\om_\zeta=\om_1^{\zeta_1}\cdots\om_n^{\zeta_n}, \qquad
\om_\zeta'=(\om_1')^{\zeta_1}\cdots(\om_n')^{\zeta_n}.\leqno(6)
$$
The following $Q$-graded structure on $U_{r,s}(\mathfrak g)$ is
necessary to develop to its weight representation theory discussed
in the sequel.

\begin{coro} [{[BGH, Corollary 2.7]}] \ {\it For any
$\zeta=\sum_{i=1}^n\zeta_i\al_i\in Q$, the defining relations
$(X2)$ in $U_{r,s}(\mathfrak g)$ take the form below:
\begin{gather*}
\om_{\zeta}\,e_i\,\om_{\zeta}^{-1}=\lg \om_i',\om_\zeta\rg\,e_i,
\qquad
\om_{\zeta}\,f_i\,\om_{\zeta}^{-1}=\lg \om_i',\om_\zeta\rg^{-1}f_i,\\
{\om_{\zeta}'}\,e_i\,{\om_{\zeta}'}^{-1}=\lg \om_\zeta',
\om_i\rg^{-1} e_i,\qquad\quad
\om_{\zeta}'\,f_i\,{\om_{\zeta}'}^{-1}=\lg
\om_\zeta',\om_i\rg\,f_i.
\end{gather*}
Then $U_{r,s}(\mathfrak n^\pm)=\bigoplus_{\eta\in
Q^+}U_{r,s}^{\pm\eta}(\mathfrak n^\pm)$ is $Q^\pm$-graded, where
$$
U_{r,s}^{\eta}(\mathfrak n^\pm)=\left\{\,a\in U_{r,s}(\mathfrak
n^\pm)\;\left|\; \om_\zeta\,a\,\om_\zeta^{-1}=\lg
\om_\eta',\om_\zeta\rg\,a, \ \om_\zeta'\,a\,{\om_\zeta'}^{-1}=\lg
\om_\zeta',\om_\eta\rg^{-1} \,a\,\right\}\right.,
$$
for $\eta\in Q^+\cup Q^-$.

\medskip
Moreover, $U=\bigoplus_{\eta\in Q}U_{r,s}^\eta(\mathfrak g)$ is
$Q$-graded such that
\begin{equation*}\begin{split}
U_{r,s}^\eta(\mathfrak g)&=\left\{\,\left.\sum
F_\al\om_{\mu}'\om_\nu E_\be\in U\; \right|\;
\om_\zeta\,(F_\al\om_{\mu}'\om_\nu E_\be)\,\om_{\zeta}^{-1}=
\lg \om'_{\beta-\al},\om_\zeta\rg\,F_\al\om_\mu'\om_\nu E_\be,\right.\\
&\quad \left.\om_\zeta'\,(F_\al\om_\mu'\om_\nu
E_\be)\,{\om_{\zeta}'}^{-1}= \lg
\om_\zeta',\om_{\beta-\al}\rg^{-1}\,F_\al\om_\mu'\om_\nu E_\be,\;
\text{\it with } \; \beta-\al=\eta\right\},
\end{split}
\end{equation*}
where $F_\al$ $($resp. $E_\be$$)$ is a certain monomial
$f_{i_1}{\cdots} f_{i_l}$ $($resp. $e_{j_1}{\cdots} e_{j_m}$$)$
such that $\al_{i_1}+{\cdots}+\al_{i_l}=\al$ $($resp.
$\al_{j_1}+{\cdots}+\al_{j_m}=\be$$)$.\hfill\qed}
\end{coro}

\bigskip
\section{Finite-Dimensional Weight Representation
Theory and Category $\mathcal O$}
\medskip

As we know, the standard triangular decomposition of
$U_{r,s}(\mathfrak g)$ suggests that $U_{r,s}(\mathfrak g)$
possesses highest weight representation theory. Indeed, this has
been developed by Benkart and Witherspoon in [BW2] for $\mathfrak
g=\mathfrak {gl}_n$ or $\mathfrak {sl}_n$. In principle, one can
expect the same theory to be valid as well for $\mathfrak
g=\mathfrak{so}_{2n+1}$, $\mathfrak {so}_{2n}$ and $\mathfrak
{sp}_{2n}$. To establish this, we will follow Benkart and
Witherspoon's main ideas. However, to treat these cases in a unified fashion, we need
to have better insights here and there in order to generalize
the techniques used in the type $A$ case. Throughout the article, we assume
that $\mathbb K$ is an algebraically closed field containing
$\mathbb Q(r,s)$ as a subfield and $rs^{-1}$ is not a root of unity.

Let $\Lambda$ be the weight lattice of $\mathfrak g$ for $\mathfrak
g=\mathfrak {so}_{2n+1}$, $\mathfrak {so}_{2n}$, or $\mathfrak
{sp}_{2n}$, respectively. We adopt similar notions and notations in
[BW1]. Associated to any $\lambda\in \Lambda$ is an algebra
homomorphism $\hat \lambda$ from the subalgebra $U^0$ over $\mathbb
K$ generated by the elements $\om_i$, $\om_i'$ ($1\le i\le n$) to
$\mathbb K$ given by
$$
\hat \lambda(\om_i)=\lg \om_\lambda',\om_i\rg, \qquad \hat
\lambda(\om_i')=\lg \om_i',\om_\lambda\rg^{-1}, \leqno(1)
$$
here we extend the definition of $\lg\, ,\rg$ from $\lambda\in Q$ to
$\lambda\in\Lambda$ via taking appropriate half-integer powers when
necessary, observing that
$\Lambda\subseteq\bigoplus_{i=1}^n\frac1{2}\mathbb
Z\al_i\subseteq\bigoplus_{i=1}^n\frac1{2}\mathbb Z\ep_i$.

Let $M$ be a $U$-module of dimension $d<\infty$ where
$U=U_{r,s}(\mathfrak g)$. As $\mathbb K$ is algebraically closed,
by linear algebra, we have
$$
M=\bigoplus_\chi M_\chi, $$ where each $\chi:\, U^0\lra \mathbb K$
is an algebra homomorphism, and $M_\chi$ is the generalized
eigenspace given by
$$
M_\chi=\left\{\,m\in M\mid
(\om_i-\chi(\om_i)1)^dm=0=(\om_i'-\chi(\om_i')1)^dm, \ \forall\;
i\,\right\}.\leqno(2)
$$
When $M_\chi\ne0$ we say that $\chi$ is a weight and $M_\chi$ is the
corresponding weight space. In the case when $M$ decomposes into
genuine eigenspaces relative to $U^0$, we say that $U^0$ acts
semisimply on $M$.

Relations in $(X2)$ imply
$$
e_jM_\chi\subseteq M_{\chi\cdot\widehat {\al_j}}\,,\qquad
f_jM_\chi\subseteq M_{\chi\cdot (\widehat{-\al_j})}\,,\leqno(3)
$$
where $\widehat{\al_j}$ is as in (1), and $\chi\cdot\psi$ is the
homomorphism with values $(\chi\cdot
\psi)(\om_i)=\chi(\om_i)\psi(\om_i)$ and $(\chi\cdot
\psi)(\om_i')=\chi(\om_i')\psi(\om_i')$. In fact, if
$(\om_i-\chi(\om_i)1)^km=0$, then $(\om_i-\chi(\om_i)\lg
\om_j',\om_i\rg 1)^ke_jm$ $=0$, and similarly for $\om_i'$ and for
$f_j$. On the one hand, (3) means that the sum of the eigenspaces
is a submodule of $M$, and so if $M$ is simple, the sum must be
$M$ itself, meanwhile we may replace the power $d$ in (2) by $1$,
that is, $U^0$ acts semisimply on each simple $M$. On the other
hand, a direct consequence of (3) is that for each simple $M$
there is a homomorphism $\chi$ so that all the weights of $M$ are
of the form $\chi\cdot \hat\zeta$, where $\zeta\in Q$.

When all the weights of a module $M$ are of the form $\hat
\lambda$, where $\lambda\in\Lambda$, we say that $M$ has weights
in $\Lambda$. Any simple $U$-module having one weight in $\Lambda$
has all its weights in $\Lambda$.

The observation below, which arises from Benkart and Witherspoon
[BW2, Proposition 3.5] in the case when $\mathfrak g=\mathfrak
{gl}_n$, or $\mathfrak {sl}_n$, also holds in our cases when
$\mathfrak g=\mathfrak{so}_{2n+1}$, $\mathfrak {so}_{2n}$ and
$\mathfrak {sp}_{2n}$.

\begin{lemma} \ For $\mathfrak g=\mathfrak {sl}_n$, $\mathfrak{so}_{2n+1}$,
$\mathfrak {so}_{2n}$ and $\mathfrak {sp}_{2n}$, suppose that
$\hat \zeta=\hat \eta$, where $\zeta,\,\eta\in \Lambda$. Assume
that $rs^{-1}$ is not a root of unity, then $\zeta=\eta$.
\end{lemma}
\begin{proof} \ The proof for $\mathfrak g=\mathfrak{sl}_n$ was given in
[BW1, Proposition 3.5]. We now give the proof case by case for
$\mathfrak g=\mathfrak {so}_{2n+1}$, $\mathfrak {sp}_{2n}$ and
$\mathfrak {so}_{2n}$, respectively.

For $\zeta=\sum_{i=1}^n\zeta_i\al_i\in\Lambda$, by definition, we
have
\begin{gather*} \hat
\zeta(\om_i)=\lg\om_\zeta',\om_i\rg=\begin{cases}
r^{2(\ep_i,\zeta)}s^{2(\ep_{i+1},\zeta)}, \quad\qquad\qquad & \quad\;\, i<n,\\
r^{2(\ep_n,\zeta)}(rs)^{-\zeta_n}. \quad\qquad\qquad & \quad\;\,
i=n;
\end{cases}\tag{\text{B}}\\
\hat \zeta(\om_i')=\lg\om_i',\om_\zeta\rg^{-1}=\begin{cases}
r^{2(\ep_{i+1},\zeta)}s^{2(\ep_i,\zeta)}, \quad\qquad\quad & \quad\;\, i<n,\\
s^{2(\ep_n,\zeta)}(rs)^{-\zeta_n}. \quad\qquad\quad & \quad\;\,
i=n.
\end{cases}\\
\hat \zeta(\om_i)=\lg\om_\zeta',\om_i\rg=\begin{cases}
r^{(\ep_i,\zeta)}s^{(\ep_{i+1},\zeta)}, \quad\qquad\qquad & \quad\;\, i<n,\\
r^{2(\ep_n,\zeta)}(rs)^{-2\zeta_n}. \quad\qquad\qquad & \quad\;\,
i=n;
\end{cases}\tag{\text{C}}\\
\hat \zeta(\om_i')=\lg\om_i',\om_\zeta\rg^{-1}=\begin{cases}
r^{(\ep_{i+1},\zeta)}s^{(\ep_i,\zeta)}, \quad\qquad\quad & \quad\;\, i<n,\\
s^{2(\ep_n,\zeta)}(rs)^{-2\zeta_n}. \quad\qquad\quad & \quad\;\,
i=n.
\end{cases}\\
\hat \zeta(\om_i)=\lg\om_\zeta',\om_i\rg=\begin{cases}
r^{(\ep_i,\zeta)}s^{(\ep_{i+1},\zeta)}, \ & \quad\;\, i<n,\\
r^{(\ep_{n-1},\zeta)}s^{-(\ep_n,\zeta)}(rs)^{-2\zeta_{n-1}}. \ &
\quad\;\, i=n;
\end{cases}\tag{\text{D}}\\
\hat \zeta(\om_i')=\lg\om_i',\om_\zeta\rg^{-1}=\begin{cases}
r^{(\ep_{i+1},\zeta)}s^{(\ep_i,\zeta)}, \ & \ i<n,\\
r^{-(\ep_n,\zeta)}s^{(\ep_{n-1},\zeta)}(rs)^{-2\zeta_{n-1}}. \ & \
i=n.
\end{cases}
\end{gather*}

Denote $\mu=\zeta-\eta$, from $\hat\zeta(\om_n)=\hat\eta(\om_n)$ and
$\hat \zeta(\om_n')=\hat\eta(\om_n')$, in the type $B$ or $C$ case,
we get $r^{2(\ep_n,\mu)}(rs)^{-\mu_n}=1$,
$s^{2(\ep_n,\mu)}(rs)^{-\mu_n}=1$; or
$r^{2(\ep_n,\mu)}(rs)^{-2\mu_n}=1$,
$s^{2(\ep_n,\mu)}(rs)^{-2\mu_n}=1$. So $(rs^{-1})^{2(\ep_n,\mu)}=1$
which, together with the assumption, means the integer
$2(\ep_n,\mu)=0$, that is,
$$
\mu_{n-1}=\mu_n, \quad (\text{\it for type \text{\rm B}}), \quad
or \qquad \mu_{n-1}=2\mu_n, \quad (\text{\it for type \text{\rm
C}}). \leqno(4)
$$
Again from $\hat\zeta(\om_{n-1})=\hat\eta(\om_{n-1})$ and $\hat
\zeta(\om_{n-1}')=\hat\eta(\om_{n-1}')$, in the type $B$ or $C$
case, we get $(\al_{n-1},\mu)=0$, that is,
$$
\mu_{n-2}=\mu_n,\quad(\text{\it for type \text{\rm B}}), \quad or
\qquad \mu_{n-2}=2\mu_n,\quad(\text{\it for type \text{\rm
C}}).\leqno(5)
$$
But similar to the deduction in the case of type $A$ (see [BW1]),
noting $\mu_0=0$, we have
\begin{gather*}
\mu_{i+2}-\mu_{i+1}-\mu_i+\mu_{i-1}=0, \quad
(i=1,2,\cdots,n-2), \tag{6}\\
\mu_{2k}=k\mu_2,\qquad \mu_{2k+1}=k\mu_2+\mu_1.\tag{7}
\end{gather*}

Thus, by (4), (5) \& (6), we get
$\mu_n=\mu_{n-1}=\cdots=\mu_1=\mu_0=0$ in the type $B$ case. For the
type $C$ case, if $n=2m$, by (5) \& (7), we get
$\mu_{n-2}=(m-1)\mu_2=2\mu_n=2m\mu_2$, i.e., $\mu_2=0$, so
$\mu_n=0$; if $n-1=2m$, then by (4), (5), \& (7), we get
$m\mu_2=\mu_{n-1}=\mu_{n-2}=(m-1)\mu_2+\mu_1$, i.e., $\mu_2=\mu_1$,
again by (4) \& (7), we get $\mu_2=0$, so $\mu_n=0$, which is
reduced to the precondition of the proof in the type $A$ case.
Hence, using the same argument as in the case of type $A$ ([BW1]),
we have $\mu=0$. Therefore, $\zeta=\eta$ in both cases $B$ and $C$.

For the type $D$ case, from $\hat \zeta(\om_i)=\hat\eta(\om_i)$ and
$\hat \zeta(\om_i')=\hat\eta(\om_i')$ for $i=n-1,\,n$, we have
$(rs^{-1})^{(\al_{n-1},\mu)}=1$ and $(rs^{-1})^{(\al_n,\mu)}=1$,
that means, together with the assumption, the integers
$(\al_{n-1},\mu)=0$ and $(\al_n,\mu)=0$. So we get
$\mu_{n-2}=2\mu_{n-1}=2\mu_n$. If $n=2m$, then
$(m-1)\mu_2=\mu_{n-2}=2m\mu_2$, i.e., $\mu_2=0$. If $n-1=2m$,
applying (7) to $\mu_{n-1}=\mu_n$, we get $\mu_1=0$; applying (7) to
$\mu_{n-2}=2\mu_{n-1}$, we get $\mu_2=0$. So we have $\mu_n=0$ for
any $n$. Using the same proof as in the case of type $A$, we obtain
$\mu=0$, i.e., $\zeta=\eta$.
\end{proof}

\begin{remark} \ Lemma 2.1 indicates that under the assumption
that $rs^{-1}$ is not a root of unity, we may simplify the notation
by writing $M_\lambda$ for the weight space rather than writing
$M_{\hat \lambda}$ for $\lambda\in\Lambda$. So it makes sense to let
(3) take the classical form: $e_jM_\lambda\subseteq
M_{\lambda+\al_j}$ and $f_jM_{\lambda}\subseteq M_{\lambda-\al_j}$.
\end{remark}

\medskip
Similar to the proof of [BW2, Corollary 3.14], we have

\begin{coro} \ {\it Let $M$ be a finite-dimensional
$U_{r,s}(\mathfrak g)$-module for $\mathfrak
g=\mathfrak{sl}_{n{+}1}$,
$\mathfrak{so}_{2n{+}1},\,\mathfrak{so}_{2n}$ or $\mathfrak
{sp}_{2n}$. Assume that $rs^{-1}$ is not a root of unity, then the
elements $e_i,\,f_i$ $(1\le i\le n)$ act nilpotently on
$M$.\hfill\qed}
\end{coro}

Obviously, when $rs^{-1}$ is not a root of unity, a
finite-dimensional simple $U$-module is a highest weight module by
Corollary 2.3 and (3).

\smallskip
We state the definition of the category $\mathcal O$ of weight
$U$-modules as in [BW1, Section 4].

\begin{definition} \ {\it Let $\mathcal O$ denote the category of
modules $M$ for $U_{r,s}(\mathfrak g)$ $($where $\mathfrak
g=\mathfrak{so}_{2n+1}$, $\mathfrak {so}_{2n}$, or $\mathfrak
{sp}_{2n})$ which satisfy the following conditions:

${(\mathcal O1)}$ \ $U^0$ acts semisimply on $M$, and the set
$\text{\rm wt}(M)$ of weights of $M$ belongs to $\Lambda:
\,M=\bigoplus_{\lambda\in\text{\rm wt}(M)}M_\lambda$, where
$M_\lambda=\{\,m\in M\mid \om_i.m=\lg \om_\lambda',\om_i\rg\, m,\;
\om_i'.m=\lg \om_i',\om_\lambda\rg^{-1} m, \;\forall \;i\ \};$

${(\mathcal O2)}$ \ $\dim_\mathbb KM_\lambda<\infty$ for all
$\lambda\in \text{\rm wt}(M);$

${(\mathcal O3)}$ \ $\text{\rm wt}(M)\subseteq\cup_{\mu\in
F}(\mu-Q^+)$ for some finite set $F\subset \Lambda$.

\noindent The morphisms in $\mathcal O$ are $U$-module
homomorphisms.}
\end{definition}

Actually, the category $\mathcal O$ just focuses on the class of
the so-called {\it type $1$} $U$-modules like in the case of
Drinfel'd-Jimbo quantum groups (see [J], [Jo], [KS]), which is
closed under taking sub-object or sub-quotient object, making
finite direct sum and taking tensor product.

\smallskip

Let $V^\psi$ be the one-dimensional $\mathcal B$-module on which
$e_i$ acts as multiplication by $0$ $(1\le i\le n)$, and $U^0$
acts via $\psi$, an algebra homomorphism from $U^0$ to $\mathbb
K$. As usual, we can define the Verma module $M(\psi)$ with
highest weight $\psi$ to be the $U$-module induced from $V^\psi$,
that is,
$$M(\psi)=U\otimes_{\mathcal B} V^{\psi}.$$
Set $v_\psi=1\otimes v\in M(\psi)$, where $v(\ne0)\in V^\psi$.
Then $e_i.v_\psi=0$ $(1\le i\le n)$ and $a.v_\psi=\psi(a)\,v_\psi$
for any $a\in U^0$ by construction. By Corollary 1.8,
$M(\psi)\cong U_{r,s}(\mathfrak n^-)\otimes v_\psi$. Corollary 1.9
indicates that each Verma module $M(\psi)\in \text{Ob}(\mathcal
O)$ if and only if $\psi\in\hat\Lambda$.

Let $N'$ be a proper submodule of $M(\psi)$, then (3) implies that
$$
N'\subset \sum_{\mu\in Q^+-\{0\}}M(\psi)_{\psi\cdot{(\widehat
{-\mu})}},
$$
as $M(\psi)_\psi=\mathbb Kv_\psi$ generates $M(\psi)$. Hence,
$M(\psi)$ has a unique maximal submodule $N$, namely the sum of all
proper submodules, and a unique simple quotient, $L(\psi)$.
Actually, all finite-dimensional simple $U$-modules are of this
form, as the Theorem below indicates (which was proved by Benkart
and Witherspoon [BW2, Theorem 2.1] in the case when $\mathfrak
g=\mathfrak {gl}_n,\,\mathfrak {sl}_n$, but still holds with the
same proof for our cases of $\mathfrak g$).

\begin{theorem}  \ For $\mathfrak g=\mathfrak
{sl}_{n+1},\,\mathfrak {so}_{2n+1},\,\mathfrak {so}_{2n}$ or
$\mathfrak {sp}_{2n}$, let $M$ be a $U_{r,s}(\mathfrak g)$-module,
on which $U^0$ acts semisimply and which contains an element $m\in
M_\psi$ $(\psi\in \text{\rm Hom}_{\text{\rm Alg}}(U^0, \mathbb
K))$ such that $e_i.m=0$ for all $i$. Then there is a unique
homomorphism of $U_{r,s}(\mathfrak g)$-modules $F:\,M(\psi)\lra M$
with $F(v_\psi)=m$. In particular, if $rs^{-1}$ is not a root of
unity and $M$ is a finite-dimensional simple $U_{r,s}(\mathfrak
g)$-module, then $M\cong L(\psi)$ for some weight
$\psi$.\hfill\qed
\end{theorem}

As in [BW2, Lemma 2.3],  it is easy to verify the commutation
relations below.

\begin{lemma} \ For \,$m\ge 1$, set
$[m]_i=\frac{r_i^m-s_i^m}{r_i-s_i}$. Then for $1\le i\le n$, we
have
\begin{gather*}
e_if_i^m=f_i^me_i+[m]_i\,f_i^{m-1}\frac{r_i^{1-m}\om_i-s_i^{1-m}\om_i'}{r_i-s_i},\\
e_i^mf_i=f_ie_i^m+[m]_i\,e_i^{m-1}\frac{s_i^{1-m}\om_i-r_i^{1-m}\om_i'}{r_i-s_i}.
\end{gather*}\hfill\qed
\end{lemma}

Set $\al^\vee=\frac{2\al}{(\al,\al)}$, for any simple root
$\al\in\Pi$, then for any $\lambda\in\Lambda$,
$(\lambda,\al^\vee)\in\mathbb Z$ by definition. Let
$\Lambda^+\subset \Lambda$ be the subset of dominant weights, that
is, $\Lambda^+=\{\,\lambda\in\Lambda\mid (\lambda,\al_i^\vee)\ge
0, \; \text{\rm for } \;1\le i\le n\,\}$.

Similar to [BW2, Lemma 2.4] in the type $A$ case, we have

\begin{lemma} \ For $\mathfrak g=\mathfrak {so}_{2n+1},\,\mathfrak {so}_{2n}$
and $\mathfrak {sp}_{2n}$, assume that $rs^{-1}$ is not a root of
unity. Let $M$ be a nonzero finite-dimensional $U_{r,s}(\mathfrak
g)$-module on which $U^0$ acts semisimply. Suppose there is some
nonzero vector $v\in M_\lambda$ with $\lambda\in\Lambda$ such that
$e_i.v=0$ for all $i$ $(1\le i\le n)$. Then $\lambda\in\Lambda^+$.
\end{lemma}
\begin{proof} \ It suffices to prove that $(\lambda,\al_n^\vee)\ge
0$, as the proof of $(\lambda,\al_i^\vee)\ge0$ ($1\le i<n$) is the
same as that of [BW2, Lemma 2.4].

Since $f_n$ acts nilpotently on $M$ by Corollary 2.3, there is
some integer $m\ge 0$ such that $f_n^{m+1}.v=0$ but $f_n^m.v\ne
0$. Applying $e_n$ to $f_n^{m+1}.v=0$, using Lemma 2.6 and the
fact that $e_n.v=0$, we get
$r_n^{-m}\hat\lambda(\om_n)=s_n^{-m}\hat\lambda(\om_n')$.
Equivalently,
\begin{gather*}
r_n^{-m}r^{2(\ep_n,\lambda)}(rs)^{-\lambda_n}=s_n^{-m}s^{2(\ep_n,\lambda)}(rs)^{-\lambda_n},\qquad\qquad\quad
\qquad\qquad\qquad
\text{({\it for type} B)} \\
r_n^{-m}r^{2(\ep_n,\lambda)}(rs)^{-2\lambda_n}=s_n^{-m}s^{2(\ep_n,\lambda)}(rs)^{-2\lambda_n},\qquad\qquad\quad
\qquad\qquad\quad\,
\text{({\it for type}  C)} \\
r^{-m}r^{(\ep_{n-1},\lambda)}s^{-(\ep_n,\lambda)}(rs)^{-2\lambda_{n-1}}=
s^{-m}r^{-(\ep_n,\lambda)}s^{(\ep_{n-1},\lambda)}(rs)^{-2\lambda_{n-1}},\quad
\text{({\it for type} D)}
\end{gather*}
or equivalently,
$$
(r_ns_n^{-1})^{-m+(\lambda,\al_n^\vee)}=1, \quad\text{({\it for
types} B, C, D)}.
$$
The assumption of $rs^{-1}$ forces $(\lambda,\al_n^\vee)=m\ge 0$.
Therefore, $\lambda\in\Lambda^+$.
\end{proof}

\begin{coro} \ {\it For $\mathfrak g=\mathfrak
{so}_{2n+1},\,\mathfrak {so}_{2n}$ and $\mathfrak {sp}_{2n}$,
assume that $rs^{-1}$ is not a root of unity, then any
finite-dimensional simple $U_{r,s}(\mathfrak g)$-module with
weights in $\Lambda$ is isomorphic to $L(\lambda)$ for some
$\lambda\in\Lambda^+$.\hfill\qed}
\end{coro}

The representation theory of $U_{r,s}(\mathfrak {sl}_2)$, developed
by Benkart and Witherspoon in [BW2], plays a crucial role in the
classification of finite-dimensional simple modules for
$U_{r,s}(\mathfrak{sl}_n)$ (see [BW2, Section 2]) like in the
classical case of the simple Lie algebras or in the quantized case
of the Drinfel'd-Jimbo quantum groups. Note the observation arising
from the structure constants of $U_{r,s}(\mathfrak g)$ for
$\mathfrak g=\mathfrak {so}_{2n+1},\,\mathfrak {so}_{2n}$ and
$\mathfrak {sp}_{2n}$: for any vertex $i$ from the corresponding
Dynkin diagram of type $B,\, C$, or $D$, respectively, $
\lg\om_i',\om_i\rg=r_is_i^{-1}$ always holds.  This fact guarantees
that even in the two-parameter quantum orthogonal or symplectic
groups $U_{r,s}(\mathfrak g)$, there exist isomorphic copies of
$U_{r,s}(\mathfrak{sl}_2)$ as well. This suggests that these quantum
groups possess a familiar finite-dimensional (weight) representation
theory provided that $rs^{-1}$ is not a root of unity.

Now let us recall the representation theory for
$U_{r,s}(\mathfrak{sl}_2)$. The first two assertions of the
following Proposition comes from [BW2, Proposition 2.8 (i)], the
last one may be regarded as an intrinsic generalization of [BW2,
Proposition 2.8 (ii)] with a deep insight.

\begin{prop}  \ Assume that $rs^{-1}$ is not a root of unity.
For $U=U_{r,s}(\mathfrak {sl}_2)$ generated by $e$, $f$, $\om$ and
$\om'$, for a given $\phi\in\text{\rm Hom}_{\text{\rm Alg}}(U^0,
\mathbb K)$, set $\phi=\phi(\om)$, $\phi'=\phi(\om')$, and in the
Verma module $M(\phi)$, put $v_j=f^j/[j]!\otimes v_\phi$ for $j\ge
0$. Then

$(\text{\rm i})$ \ $M(\phi)$ is a simple $U$-module if and only if
$\phi\cdot r^{-j}-\phi'\cdot s^{-j}\ne 0$ for any $j\ge 0$.

$(\text{\rm ii})$ \ If $\phi(\om')=\phi(\om)(rs^{-1})^{-m}$ for
some integer $m\ge 0$, then $\text{\rm Span}_\mathbb K\{\,v_j\mid
j\ge m+1\,\}\cong M(\phi-(m+1)\al)$ is the unique maximal
submodule of $M(\phi)$. The quotient is the $(m+1)$-dimensional
simple module $L(\phi)$ spanned by vectors $v_0,v_1,\cdots,v_m$
and having $U$-action given by
\begin{gather*}
\om.v_j=\phi\cdot(rs^{-1})^{-j}v_j,\qquad
\om'.v_j=\phi\cdot(rs^{-1})^{-(m-j)}v_j,\\
e.v_j=\phi\cdot r^{-m}[m+1-j]\,v_{j-1}, \quad(v_{-1}=0)\tag{8}\\
f.v_j=[j+1]\,v_{j+1}. \quad (v_{m+1}=0)
\end{gather*}
Any $(m+1)$-dimensional simple $U$-module is isomorphic to
$L(\phi)$ for some such $\phi$.

$(\text{\rm iii})$ \ If
$\nu=\nu_1\lambda_1+\cdots+\nu_n\lambda_n\in\Lambda^+$, where
$\lambda_i$ is the $i$-th fundamental weight for $\mathfrak g$,
then $\hat \nu(\om_i')=\hat\nu(\om_i)(r_is_i^{-1})^{-\nu_i}$, and
the $U_i$-module $L(\nu_i\lambda_i)$ is $(\nu_i+1)$-dimensional
and has $U_i$-action given by $(8)$ with $\phi_i=\hat\nu(\om_i)$,
where $U_i$ is the copy of $U_{r,s}(\mathfrak {sl}_2)$ in
$U_{r,s}(\mathfrak g)$ corresponding to the $i$-th vertex of the
Dynkin diagram.
\end{prop}
\begin{proof} \ For the proof of the last assertion, it suffices to
show that there hold
$$
\frac{\hat\nu(\om_i')}{\hat\nu(\om_i)}=
(r_is_i^{-1})^{-(\al_i^\vee,\,\nu)}=\frac{\widehat{\nu_i\lambda_i}(\om_i')}{\widehat{\nu_i\lambda_i}(\om_i)},
\quad (\text{\it for any } \; i)\leqno(9)
$$
for $\mathfrak g=\mathfrak{sl}_n,\, \mathfrak
{so}_{2n+1},\,\mathfrak {so}_{2n}$ and $\mathfrak {sp}_{2n}$.

In the type $A$ case, we have $\lambda_i=\ep_1+\cdots+\ep_i$ for
$1\le i\le n$ and $\lambda_n=0$. By definition,
\begin{equation*}
\begin{split}
\frac{\hat\nu(\om_i')}{\hat\nu(\om_i)}&=\frac{s^{(\ep_i,\nu)}r^{(\ep_{i+1},\nu)}}{r^{(\ep_i,\nu)}s^{(\ep_{i+1},\nu)}}
=(rs^{-1})^{-(\al_i,\nu)}=(rs^{-1})^{-\nu_i}\\
&=\frac{s^{(\ep_i,\nu_i\lambda_i)}r^{(\ep_{i+1},\nu_i\lambda_i)}}{r^{(\ep_i,\nu_i\lambda_i)}s^{(\ep_{i+1},\nu_i\lambda_i)}}
=\frac{\widehat{\nu_i\lambda_i}(\om_i')}{\widehat{\nu_i\lambda_i}(\om_i)}.
\end{split}
\end{equation*}

For types $B,\,C$ and $D$, it suffices to consider types $B_2$,
$C_2$ and $D_4$, respectively.

In the type $B_2$ case, we have $\lambda_1=\ep_1$,
$\lambda_2=\frac1{2}(\ep_1+\ep_2)$. By the defining formula (B) in
Lemma 2.1, for $i=1$, it follows directly from the argument in the
type $A$ case; while for $i=2$, we get
$$
\frac{\hat\nu(\om_2')}{\hat\nu(\om_2)}=\frac{s^{2(\ep_2,\nu)}}{r^{2(\ep_2,\nu)}}=(rs^{-1})^{-(\al_2^\vee,\,\nu)}=
\frac{s^{2(\ep_2,\nu_2\lambda_2)}}{r^{2(\ep_2,\nu_2\lambda_2)}}
=\frac{\widehat{\nu_2\lambda_2}(\om_2')}{\widehat{\nu_2\lambda_2}(\om_2)}.
$$

In the type $C_2$ case, we have $\lambda_1=\ep_1$,
$\lambda_2=\ep_1+\ep_2$. It suffices to consider the case $i=2$.
Similarly, we have
$$
\frac{\hat\nu(\om_2')}{\hat\nu(\om_2)}=\frac{s^{2(\ep_2,\nu)}}{r^{2(\ep_2,\nu)}}=(r_2s_2^{-1})^{-(\al_2^\vee,\,\nu)}=
\frac{s^{2(\ep_2,\nu_2\lambda_2)}}{r^{2(\ep_2,\nu_2\lambda_2)}}
=\frac{\widehat{\nu_2\lambda_2}(\om_2')}{\widehat{\nu_2\lambda_2}(\om_2)}.
$$

In the type $D_4$ case,  we have
$\lambda_1=\ep_1,\,\lambda_2=\ep_1+\ep_2,\,\lambda_3=\frac1{2}(\ep_1+\ep_2+\ep_3-\ep_4),\,
\lambda_4=\frac1{2}(\ep_1+\ep_2+\ep_3+\ep_4)$.  It suffices to
consider the cases $i=3,\,4$. By the formula (D) in Lemma 2.1, we
have
$$
\frac{\hat\nu(\om_i')}{\hat\nu(\om_i)}=(rs^{-1})^{(\al_i,\,\nu)}=\frac{\widehat{\nu_i\lambda_i}(\om_i')}{\widehat{\nu_i\lambda_i}(\om_i)},
$$
for $i=3,\,4$.

The proof is completed.
\end{proof}

Proposition 2.9 (iii) and its proof imply the following result.

\begin{coro} \ {\it Assume that $rs^{-1}$ is not a root of
unity and $\lambda\in\Lambda^+$, set $\nu_i=(\lambda,
\al_i^\vee)$, then each vector $f_i^{\nu_i+1}.v_\lambda$ in the
Verma $U$-module $M(\lambda)$ generates the Verma submodule
$M(\lambda-(\nu_i+1)\al_i)$ for all $i$, where $\mathfrak
g=\mathfrak{sl}_{n+1},\, \mathfrak {so}_{2n+1},\,\mathfrak
{so}_{2n}$ or $\mathfrak {sp}_{2n}$.}
\end{coro}
\begin{proof} \ It follows from a direct calculation of
$e_if^{\nu_i+1}_i.v_\lambda=0$ by Lemma 2.6 and (9).
\end{proof}

More generally, we have

\begin{prop} \ Let $M(\lambda)$ be a Verma module with
$\lambda\in\Lambda^+$. Then for every element $\omega$ of the Weyl
group $\mathcal W$ of $\mathfrak g$, there exists a Verma
submodule in $M(\lambda)$ with highest weight
$$
\lambda_\omega=\omega(\lambda+\rho)-\rho,\leqno(10)
$$
where $\rho$ is the half-sum of all positive roots of $\mathfrak
g$. Every simple $U$-module as a composition factor of
$M(\lambda)$ determines a highest weight module in $\mathcal O$.
These highest weights are of the form $(10)$.
\end{prop}
\begin{proof} \ The proof of this proposition is analogous to that
of the corresponding assertion in the classical theory (see
Dixmier [D]).
\end{proof}

\begin{lemma} \ For any simple $U$-module $L(\lambda)$ with
$\lambda\in\Lambda^+$, take any $\be=\sum_{i=1}^nm_i\al_i\in Q^+$
such that $m_i\le (\lambda,\al_i^\vee)$, $\forall\; i$, then the
linear mapping $U^{-\be}_{r,s}(\mathfrak n^-)\ni x\mapsto
x.v_\lambda$ is injective.
\end{lemma}
\begin{proof} \ By the definition of the Verma module, it is enough
to show that $\lambda-\be$ is not a weight of the maximal
$U$-submodule $N$. This follows from Proposition 2.11, because no
set of weights $\{\,\lambda_\omega-\sum_{i=1}^nn_i\al_i\mid
n_i\in\mathbb Z^+\,\}$, $\omega\in {\mathcal W}-\{1\}$ contains
$\lambda-\be$.
\end{proof}

\begin{lemma} \ If an element $a\in
U_{r,s}^{-\be}(\mathfrak n^-)$ satisfies the relations $e_ia=ae_i$
for $i=1,2,\cdots,n$, then we have $a=0$. If $f_ib=bf_i$,
$i=1,2,\cdots,n$, for some $b\in U^\be_{r,s}(\mathfrak n)$, then
$b=0$.
\end{lemma}
\begin{proof} \ Write $\be=\sum_{i=1}^nm_i\al_i\in Q^+$, and take a
dominant weight $\lambda\in\Lambda^+$ such that
$(\lambda,\al_i^\vee)\ge m_i$ for all $i$. Consider the simple
$U$-module $L(\lambda)$ with highest weight vector $v_\lambda$.
Since $(e_ia).v_\lambda=(ae_i).v_\lambda=0$ for all $i$, the
vector $a.v_\lambda$ generates a proper submodule of $L(\lambda)$.
Thus $a.v_\lambda=0$, as $L(\lambda)$ is simple. Hence $a=0$ by
Lemma 2.12.

In order to prove the second assertion, we introduce a $\mathbb
Q$-algebra isomorphism $\theta:\,U_{r,s}(\mathfrak g)\lra
U_{r,s}(\mathfrak g)$ defined by
\begin{gather*}
\theta(r)=s^{-1},\qquad \theta(s)=r^{-1},\\
\theta(\om_i)=\om_i',\qquad \theta(\om_i')=\om_i,\tag{11}\\
\theta(e_i)=f_i,\qquad \theta(f_i)=(r_is_i)e_i.
\end{gather*}
In fact, we can find that the image of $\theta$ is $\mathbb
Q$-algebraically isomorphic to the associated quantum group
$U_{s^{-1},r^{-1}}(\mathfrak g)$, i.e., $\text{Im}(\theta)\cong
(U_{s^{-1},r^{-1}}(\mathfrak g), \lg \,|\,\rg)$, where the pairing
$\lg\om_i'|\,\om_j\rg$ is defined via substituting $(r, s)$ by
$(s^{-1}, r^{-1})$ in the defining formula for $\lg
\om_i',\om_j\rg$ (see formulae $(1_X)$ and (2) in Section 1).

Now applying the $\mathbb Q$-algebra isomorphism $\theta$ to the
equation $f_ib=bf_i$, we get $\theta(b)=0$, by the first
assertion. Hence, $b=0$.
\end{proof}

Returning to the pairing $\lg\,,\rg:\;\mathcal B'\times \mathcal
B\lra\mathbb K$ in Proposition 1.5, and combining with the
$Q$-gradation on $U$ (see Corollary 1.9), we have

\begin{prop} \ For any $\beta\in Q^+$, the restriction of the
pairing $\lg\,,\rg$ in Proposition 1.5 to ${\mathcal
B'}^{-\be}\times \mathcal B^\be$ is nondegenerate.
\end{prop}
\begin{proof} \ We have to show that for any $a\in {\mathcal
B'}^{-\be}$ such that $\lg a, b\rg=0$ for some $b\in \mathcal
B^\be$, implies that $b=0$. This will be proved by induction with
respect to the usual ordering of $Q_+$. If $\be$ is a simple root,
then it is true by formula (2) in Section 1. Let $\be>0$ with
$\text{ht}(\be)>1$ and suppose that it holds for all $\gamma\in Q^+$
such that $\be-\gamma\in Q^+$.

Note that using the defining properties of skew-dual pairing and
the comultiplication in $U$ (see Proposition 1.2), we may check by
induction:
\begin{gather*} \lg c\,\om_\nu',\om_\mu\,d\rg=\lg
\om_\nu',\om_\mu\rg\,\lg c, d\rg, \qquad \forall \;c\in
U_{r,s}(\mathfrak n^-),\; d\in
U_{r,s}(\mathfrak n),\tag{12}\\
\lg c, d\rg=0,\qquad c\in U^{-\sigma}_{r,s}(\mathfrak n^-), \ d\in
U^{\delta}_{r,s}(\mathfrak n), \ \sigma,\,\delta\in Q^+,\
\sigma\ne\delta.\tag{13}
\end{gather*}
It suffices to assume that $b\in U^\be_{r,s}(\mathfrak n)$. By
Proposition 1.2, we can write
$$
\Delta(b)=\sum_{0\le \gamma\le\be}(\om_{\gamma}\ot 1)\,b_{\gamma},
\qquad b_{\gamma}\in U^\gamma_{r,s}(\mathfrak n)\ot
U^{\be-\gamma}_{r,s}(\mathfrak n),\leqno(14)
$$
where $b_0=b\ot 1$ and $b_\be=1\ot b$. Let $\gamma\in Q^+$,
$0<\gamma<\be$, $x\in {\mathcal B'}^{-\gamma}$ and $y\in {\mathcal
B'}^{-(\be-\gamma)}$.

By (2), (12) \& (13), we have
$$
0=\lg xy,b\rg=\lg x\ot y,\Delta(b)\rg=\lg x\ot y, (\om_{\gamma}\ot
1)\,b_{\gamma}\rg=\lg x\ot y, b_\gamma\rg.\leqno(15)
$$
By assumption, for any $\gamma'<\be$ the restriction of
$\lg\,,\rg$ to ${\mathcal B'}^{-\gamma'}\times \mathcal
B^{\gamma'}$ is nondegenerate, so is its extension to a bilinear
form on $ [\,{\mathcal B'}^{-\gamma}\ot{\mathcal
B'}^{-(\be-\gamma)}\,]\times [\,\mathcal B^\gamma\ot\mathcal
B^{\be-\gamma}\,]$.  Hence it follows from (15) that $b_\gamma=0$.
Because of (14) this means that $\Delta(b)=b\ot 1+\om_\be\ot b$.
By (13), together with $\Delta(f_i)=1\ot f_i+f_i\ot \om_i'$, we
get $\hat f_i\,\hat b=\hat b\,\hat f_i$, and then $f_i\,b=b\,f_i$
for any $i$, after using $\varphi$ (see the proof of [BGH, Theorem
2.5]). Thus, by Lemma 2.12, $b=0$.

Similar reasoning indicates that for any $b\in \mathcal B^\be$
such that $\lg a, b\rg=0$ for some $a\in{\mathcal B'}^{-\be}$
implies that $a=0$.
\end{proof}

In what follows, we consider the finite-dimensionality question of
the simple $U_{r,s}(\mathfrak g)$-modules $L(\lambda)$ with
$\lambda\in\Lambda^+$. This problem has been solved by Benkart and
Witherspoon in [BW2, Section 2] in the case when $\mathfrak
g=\mathfrak {gl}_n$, or $\mathfrak{sl}_n$. The same idea can be
used to prove that $M(\lambda)$ has a $U_{r,s}(\mathfrak
g)$-submodule $M'(\lambda)$ of finite codimension, as $L(\lambda)$
is the quotient of $M(\lambda)$ by its unique maximal submodule,
where $M'(\lambda)$ is defined by
$$
M'(\lambda)=\sum_{i=1}^nU_{r,s}(\mathfrak
g)f_i^{k_i+1}.v_\lambda\cong
\sum_{i=1}^nM(\lambda-\bigl(k_i{+}1)\al_i\bigr), \leqno(16)
$$
where $k_i=(\lambda,\al_i^\vee)$ for all $i$. That is, to prove
the module $L'(\lambda)=M(\lambda)/M'(\lambda)$ is nonzero and
finite-dimensional. $L'(\lambda)\ne0$ is clear, since any weight
in $M'(\lambda)$ is less than or equal to $\lambda-(k_i{+}1)\al_i$
for some $i$, $v_\lambda\not\in M'(\lambda)$.

\begin{lemma} $(\text{\rm i})$ $([{\rm BW2, Lemma\; 2.10}])$ \ The
elements $e_j,\,f_j$ $(1\le j\le n)$ act locally nilpotently on
$U_{r,s}(\mathfrak g)$-module $L'(\lambda)$.

$(\text{\rm ii})$ $([{\rm BW2, Lemma \;2.11}])$ \ Assume that
$rs^{-1}$ is not a root of unity, $V=\bigoplus_{j\in\mathbb
Z^+}V_{\lambda-j\al}\in\text{\rm Ob}(\mathcal O)$ is a
$U_{r,s}(\mathfrak{sl}_2)$-module for some weight
$\lambda\in\Lambda$. If $\,e,\,f$ act locally nilpotently on $V$,
then $\dim_\mathbb K V<\infty$, and the weights of $V$ are preserved
under the simple reflection taking $\al$ to $-\al$.
\end{lemma}
\begin{proof} \ The proof of (i) is parallel to the type $A$ case;
the second part assertion is direct from [BW2].
\end{proof}

\begin{prop} \ Assume that $rs^{-1}$ is not a root of unity. Then
for the $U_{r,s}(\mathfrak g)$-module $L'(\lambda)\in\text{\rm
Ob}(\mathcal O)$ with $\lambda\in\Lambda^+$, we have $\dim_\mathbb K
L'(\lambda)<\infty$, so $\dim_\mathbb K L(\lambda)$ $<\infty$.
\end{prop}
\begin{proof} \ Consider $L'(\lambda)$ as a $U_i$-module, where
$U_i$ is the copy generated by $e_i,\,f_i,\,\om_i$, $\om_i'$. For
$\mu$ a weight of $L'(\lambda)$, applying Lemma 2.15 to the
$U_i$-module
$$L'_i(\mu)=U_i.L'(\lambda)_\mu=\bigoplus_{j\in\mathbb Z^+}L'_i(\mu)_{\lambda'-j
\al_i}$$ for some weight $\lambda'\le \lambda$, we get that the
simple reflection $w_i$ preserves the weights of $L'_i(\mu)$, so
$w_i(\mu)$ is a weight of $L'(\lambda)$. That is, the Weyl group
$\mathcal W$ of $\mathfrak g$ preserves the set of weights of
$L'(\lambda)$. From Lie theory, we know that each $\mathcal W$-orbit
only contains one dominant weight. But there are only finitely many
dominant weights $\le \lambda$, and as each weight space of
$L'(\lambda)$ is of finite-dimension, we have $\dim_\mathbb
KL'(\lambda)<\infty$. Thereby, $\dim_\mathbb KL(\lambda)<\infty$.
\end{proof}

For $\mathfrak
g=\mathfrak{sl}_{n+1},\,\mathfrak{so}_{2n+1},\,\mathfrak{so}_{2n}$
or $\mathfrak {sp}_{2n}$, Corollary 2.8 and Proposition 2.16 imply
the following
\begin{coro} \ {\it A finite-dimensional
simple object in the category $\mathcal O$ is precisely a
$U_{r,s}(\mathfrak g)$-module $L(\lambda)$ for some
$\lambda\in\Lambda^+$, and $L(\lambda)\cong L(\mu)$ if and only if
$\lambda=\mu$.\hfill\qed}
\end{coro}

\subsection*{Finite-dimensional simple (weight) modules of generic type}

As noted in [BW2, Section 2], for $\mathfrak
g=\mathfrak{gl}_n,\,\mathfrak{sl}_n$, Benkart and Witherspoon gave a
description of a classification of finite-dimensional simple
$U_{r,s}(\mathfrak g)$-modules. We find that a similar structural
feature for finite-dimensional simple $U_{r,s}(\mathfrak g)$-modules
also holds when $\mathfrak g=\mathfrak
{so}_{2n+1},\,\mathfrak{so}_{2n}$, or $\mathfrak{sp}_{2n}$, after
modifying some of the treatments.

Given a one-dimensional $U_{r,s}(\mathfrak g)$-module $L$, Theorem
2.5 indicates that $L=L(\chi)$ for some $\chi\in
\text{Hom}_{\text{Alg}}(U^0,\mathbb K)$ with the elements
$e_i,\,f_i$ $(1\le i\le n)$ trivially acting on $L(\chi)$.
Relation $(X3)$ $(X=B,\,C,\,D)$ gives
$$
\chi(\om_i)=\chi(\om_i'), \qquad(1\le i\le n).\leqno(17)
$$
Conversely, if $\chi\in\text{Hom}_{\text{Alg}}(U^0,\mathbb K)$
satisfies the equation (17), then Proposition 2.9 (ii) guarantees
$\dim_{\mathbb K}L(\chi)=1$. We denote by $L_\chi$ the
one-dimensional $U_{r,s}(\mathfrak g)$-module $L(\chi)$.

The following Lemma was proved by Benkart and Witherspoon in the
case of type $A$. We will give a unified proof for the classical
types of $\mathfrak g$ based on an intrinsic observation in
Proposition 2.9 (ii) \& (iii).

\begin{lemma} \ Assume $rs^{-1}$ is not a root of unity. Given a
finite-dimensional simple $U_{r,s}(\mathfrak g)$-module $L(\psi)$
with highest weight $\psi$, there exists a pair $(\chi, \lambda)$,
where $\chi\in\text{\rm Hom}_{\text{\rm Alg}}(U^0,\mathbb K)$ such
that $(17)$ holds, and $\lambda\in\Lambda^+$, so that
$\psi=\chi\cdot\hat\lambda$, and $\text{\rm wt}\bigl(L(\psi)\bigr)
\subseteq\chi\cdot\hat\Lambda$.
\end{lemma}
\begin{proof} \ As $L(\psi)$ is finite-dimensional and simple, for
each pair of eigenvalues $(\psi(\om_i),\psi(\om_i'))$ when
considering $L(\psi)$ as a $U_i$-module (where $U_i$ is a
$U_{r,s}(\mathfrak {sl}_2)$-copy of $U_{r,s}(\mathfrak g)$),
Proposition 2.9 (ii) tells us that there exists a nonnegative
integer $\nu_i$ for each index $i$ such that
$\psi(\om_i')=\psi(\om_i)(r_is_i^{-1})^{-\nu_i}$. Set
$\lambda=\sum_{i=1}^n\nu_i\lambda_i$ where $\lambda_i$ is the $i$th
fundamental weight of $\mathfrak g$, then $\lambda\in\Lambda^+$.

Now we take $\chi(\om_i)=\psi(\om_i)\hat\lambda_i(\om_i)^{-1}$ and
$\chi(\om_i')=\psi(\om_i')\hat\lambda_i(\om_i')^{-1}$, that is,
$\chi=\psi\cdot\hat\lambda_i^{-1}\in\text{Hom}_{\text{Alg}}(U^0,
\mathbb K)$ and satisfies
\begin{equation*}
\begin{split}
\chi(\om_i')&=\psi(\om_i')\hat\lambda_i^{-1}(\om_i')=
\psi(\om_i)(r_is_i^{-1})^{-\nu_i}\hat\lambda_i^{-1}(\om_i')\\
&=\psi(\om_i)\hat\lambda^{-1}(\om_i) \qquad \bigl(\;\text{\rm by} \ (9)\;\bigr)\\
&=\chi(\om_i),
\end{split}
\end{equation*}
as required. The last assertion that $\text{\rm
wt}(L(\psi))\subseteq \chi\cdot\hat \Lambda$ is quite clear.
\end{proof}

Similar to [BW2, Theorem 2.19], for $\mathfrak
g=\mathfrak{so}_{2n+1},\,\mathfrak{so}_{2n}$, or $\mathfrak
{sp}_{2n}$, we have the classification Theorem for
finite-dimensional simple $U_{r,s}(\mathfrak g)$-modules as
follows.

\begin{theorem} \ Let $rs^{-1}$ be a non-root of unity.
Each finite-dimensional simple $U_{r,s}(\mathfrak g)$-module
$L(\psi)$ with $\psi\in\text{\rm Hom}_{\text{\rm Alg}}(U^0, \mathbb
K)$ is isomorphic to $ L_\chi\otimes L(\lambda)$, where
$\chi\in\text{\rm Hom}_{\text{\rm Alg}}(U^0, \mathbb K)$ with
$\chi(\om_i)=\chi(\om_i')$ $(1\le i\le n)$ and
$\lambda\in\Lambda^+$.\hfill\qed
\end{theorem}

\medskip
\section{$R$-matrices, Quantum Casimir Operators, Complete Reducibility}
\medskip

For any two objects $M,\,M'\in\text{Ob}(\mathcal O)$,  Benkart and
Witherspoon in [BW1, Section 4] constructed a $U_{r,s}(\mathfrak
{sl}_n)$-module isomorphism
$$R_{M',M}: \,M'\otimes M\lra M\otimes
M'$$ by a remarkable method due to Jantzen [J, Chap. 7] for the
quantum groups $U_q(\mathfrak g)$ of Drinfel'd-Jimbo type.

The aim of this section is to generalize this result to the
setting of $\mathfrak
g=\mathfrak{so}_{2n+1},\,\mathfrak{so}_{2n},\,\mathfrak
{sp}_{2n}$.

Noting that the weight lattice
$$
\Lambda\subseteq\bigoplus_{i=1}^n\frac1{2}\mathbb
Z\al_i\subseteq\bigoplus_{i=1}^n\frac1{2}\mathbb Z\ep_i,
$$
as it was done in formula (1) of Section 2, for
$\lambda\in\Lambda$, we have an algebra homomorphism $\hat
\lambda\in\text{Hom}_{\text{Alg}}(U^0, \mathbb K)$. Furthermore,
we extend the pairing $\lg\,,\,\rg$ to $\Lambda\times \Lambda$,
such that for any
$\lambda=\sum_{i=1}^np_i\al_i,\,\mu=\sum_{i=1}^nq_i\al_i\in\Lambda$
with $p_i,\,q_i\in\frac1{2}\mathbb Z$, we define
$$\lg\om_\lambda',\om_\mu\rg=\Pi_{i=1}^n\hat\lambda(\om_i)^{q_i},\leqno(1)$$
which is well-defined in the algebraically closed field $\mathbb
K$.

Now we define the map $f:\,\Lambda\times\Lambda\lra\mathbb K^*$ by
$$
f(\lambda,\mu)=\lg\om_\mu',\om_\lambda\rg^{-1}, \leqno(2)
$$
which satisfies
\begin{gather*}
f(\lambda+\mu,\nu)=f(\lambda,\nu)\,f(\mu,\nu),\\
f(\lambda,\mu+\nu)=f(\lambda,\mu)\,f(\lambda,\nu),\tag{3}\\
f(\al_i,\mu)=\lg \om_\mu',\om_i\rg^{-1}, \qquad
f(\lambda,\al_i)=\lg\om_i',\om_\lambda\rg^{-1}.
\end{gather*}
 And we define the linear
transformation $\tilde f=\tilde f_{M,M'}:\,M\otimes M'\lra
M\otimes M'$ by $\tilde f(m\otimes m')=f(\lambda,\mu)\,(m\otimes
m')$ for $m\in M_\lambda$ and $m'\in M_{\mu}'$.

Owing to $\Delta(e_i)=e_i\otimes 1+\om_i\otimes e_i$, we have
$\Delta(x)\in\sum_{0\le\nu\le\zeta}U_{r,s}^{\zeta-\nu}(\mathfrak
n)\om_\nu\otimes U_{r,s}^{\nu}(\mathfrak n)$, for all $x\in
U_{r,s}^{\zeta}(\mathfrak n)$, by induction. For each $i$, the
expression of $\Delta(x)$ defines two skew-derivations
$\partial_i,\,_i\partial: \,U_{r,s}^{\zeta}(\mathfrak n)\lra
U_{r,s}^{\zeta-\al_i}(\mathfrak n)$ such that
\begin{equation*}
\begin{split}
\Delta(x)&=x\otimes 1+\sum_{i=1}^{n}\partial_i(x)\,\om_i\otimes
e_i+\text{\it the rest},\\
\Delta(x)&=\om_\zeta\otimes
x+\sum_{i=1}^{n}e_i\,\om_{\zeta-\al_i}\otimes\,_i\partial(x)+\text{
\it the rest},
\end{split}\tag{4}
\end{equation*}
where in each case ``{\it the rest}" refers to terms involving
products of more than one $e_j$ in the second (resp. first)
factor. More precisely, parallel to [BW1, Lemma 4.6] or comparing
with [KS, Lemmas 6.14, 6.17], we have

\begin{lemma} \ For all $x\in U_{r,s}^\zeta(\mathfrak
n)$, $x'\in U_{r,s}^{\zeta'}(\mathfrak n)$, and $y\in
U_{r,s}(\mathfrak n^-)$, the following hold:

$(\text{\rm i})$ \ \;
$\partial_i(xx')=\lg\om_{\zeta'}',\om_i\rg\,\partial_i(x)\,x'+x\,\partial_i(x')$.

$(\text{\rm ii})$ \;
$_i\partial(xx')=\,_i\partial(x)\,x'+\lg\om_i',\om_\zeta\rg\,x\,_i\partial(x')$.

$(\text{\rm iii})$ \ $\lg f_iy,\, x\rg =\lg f_i, e_i\rg\,\lg y,\,
_i\partial(x)\rg=(s_i-r_i)^{-1}\lg y,\, _i\partial(x)\rg$.

$(\text{\rm iv})$ \ $\lg yf_i,\, x\rg =\lg f_i, e_i\rg\,\lg y,\,
\partial_i(x)\rg=(s_i-r_i)^{-1}\lg y,\, \partial_i(x)\rg$.

$(\text{\rm v})$ \ \,
$f_ix-xf_i=(s_i-r_i)^{-1}\bigl(\partial_i(x)\,\om_i-\om_i'\,_i\partial(x)\bigr)$.\hfill\qed
\end{lemma}

Also, for each $i$, the expression of $\Delta(y)$ for $y\in
U_{r,s}^{-\zeta}(\mathfrak n^-)$ defines two skew-derivations
$\partial_i,\,_i\partial: \,U_{r,s}^{-\zeta}(\mathfrak n^-)\lra
U_{r,s}^{-\zeta+\al_i}(\mathfrak n^-)$ such that
\begin{equation*}
\begin{split}
\Delta(y)&=y\otimes \om_\zeta'+\sum_{i=1}^n\partial_i(y)\otimes
f_i\,\om'_{\zeta-\al_i}+\text{\it
the rest},\\
\Delta(y)&=1\otimes
y+\sum_{i=1}^nf_i\otimes\,_i\partial(y)\,\om_i'+\text{\it the
rest}.
\end{split}\tag{5}
\end{equation*}
Parallel to [BW1, Lemma 4.8], we have

\begin{lemma} \ For all $y\in U_{r,s}^{-\zeta}(\mathfrak n^-)$,
$y'\in U_{r,s}^{-\zeta'}(\mathfrak n^-)$, and $x\in
U_{r,s}(\mathfrak n)$, the following hold:

$(\text{\rm i})$ \ \;
$\partial_i(yy')=\partial_i(y)\,y'+\lg\om_{\zeta}',\om_i\rg\,y\,\partial_i(y')$.

$(\text{\rm ii})$ \;
$_i\partial(yy')=\lg\om_i',\om_{\zeta'}\rg\,_i\partial(y)\,y'+y\,_i\partial(y')$.

$(\text{\rm iii})$ \ $\lg y,\,e_ix\rg =\lg f_i, e_i\rg\,\lg
\partial_i(y),\,x\rg=(s_i-r_i)^{-1}\lg \partial_i(y),\, x\rg$.

$(\text{\rm iv})$ \ $\lg y,\, xe_i\rg =\lg f_i, e_i\rg\,\lg
\,_i\partial(y),\, x\rg=(s_i-r_i)^{-1}\lg \,_i\partial(y),\,
x\rg$.

$(\text{\rm v})$ \ \,
$e_iy-ye_i=(r_i-s_i)^{-1}\bigl(\om_i\,\partial_i(y)-\,_i\partial(y)\,\om_i'\bigr)$.\hfill\qed
\end{lemma}

By Proposition 2.14, the spaces $U_{r,s}^\zeta(\mathfrak n)$ and
$U_{r,s}^{-\zeta}(\mathfrak n^-)$ are non-degenerately paired. We
may select a basis $\{u_k^\zeta\}_{k=1}^{d_\zeta}$, ($d_\zeta=\dim
U_{r,s}^{\zeta}(\mathfrak n)$), for $U_{r,s}^{\zeta}(\mathfrak n)$
and a dual basis $\{v_k^\zeta\}_{k=1}^{d_\zeta}$ for
$U_{r,s}^{-\zeta}(\mathfrak n^-)$. Then for each $x\in
U_{r,s}^{\zeta}(\mathfrak n)$ and $y\in U_{r,s}^{-\zeta}(\mathfrak
n^-)$, we have
$$
x=\sum_{k=1}^{d_\zeta}\lg v_k^\zeta, x\rg\,u_k^\zeta, \qquad
y=\sum_{k=1}^{d_\zeta}\lg y, u_k^\zeta\rg\,v_k^\zeta.\leqno(6)
$$

For $\zeta\in Q^+=\bigoplus_{i=1}^n\mathbb Z^+\al_i$, we define
$$
\Theta_\zeta=\sum_{k=1}^{d_\zeta}v_k^{\zeta}\otimes
u_k^{\zeta}.\leqno(7)
$$
Set $\Theta_\zeta=0$ if $\zeta\not\in Q^+$. Similar to [BW1, Lemma
4.10], for the cases when $\mathfrak
g=\mathfrak{so}_{2n+1},\,\mathfrak{so}_{2n}$ and
$\mathfrak{sp}_{2n}$, we also have

\begin{lemma} \ For $1\le i\le n$, the following relations hold

$(\text{\rm i})$ \quad
$(\om_i\otimes\om_i)\,\Theta_\zeta=\Theta_\zeta\,(\om_i\otimes\om_i)$,
\qquad
$(\om_i'\otimes\om_i')\,\Theta_\zeta=\Theta_\zeta\,(\om_i'\otimes\om_i');$

$(\text{\rm ii})$ \ \ $(e_i\otimes 1)\,\Theta_\zeta+(\om_i\otimes
e_i)\,\Theta_{\zeta-\al_i}=\Theta_\zeta\,(e_i\otimes
1)+\Theta_{\zeta-\al_i}\,(\om_i'\otimes e_i);$

$(\text{\rm iii})$ \ \,$(1\otimes f_i)\,\Theta_\zeta+(f_i\otimes
\om_i')\,\Theta_{\zeta-\al_i}=\Theta_\zeta\,(1\otimes
f_i)+\Theta_{\zeta-\al_i}\,(f_i\otimes \om_i).$\hfill\qed
\end{lemma}

Now we define
$$
\Theta=\sum_{\zeta\in Q^+}\Theta_\zeta.\leqno(8)
$$
Given $U_{r,s}(\mathfrak g)$-module $M$ and $M'$ in $\mathcal O$,
we apply $\Theta$ to their tensor product:
$$\Theta=\Theta_{M,M'}:\,M\otimes M'\lra M\otimes M'.$$
Note that $\Theta_\zeta: \,M_\lambda\otimes M_\mu'\lra
M_{\lambda-\zeta}\otimes M_{\mu+\zeta}'$ for all
$\lambda,\,\mu\in\Lambda$, and there are only finitely many
$\zeta\in Q^+$ such that $M_{\mu+\zeta}'\ne 0$, thanks to
condition $(\mathcal O3)$. So $\Theta$ is a well-defined linear
transformation on $M\otimes M'$. After appropriately ordering the
chosen countable bases of weight vectors for both $M$ and $M'$, we
see that each $\Theta_\zeta$ with $\zeta>0$ has a strictly
triangular matrix, while $\Theta_0=1\otimes 1$ acts as the
identity transformation on $M\otimes M'$, hence $\Theta_{M,M'}$ is
an invertible transformation.

\begin{theorem} \ Let $M$ and $M'$ be $U_{r,s}(\mathfrak
g)$-modules in $\mathcal O$ where $\mathfrak g=\mathfrak
{so}_{2n+1},\,\mathfrak{so}_{2n}$ or $\mathfrak {sp}_{2n}$. Then
the map
$$R_{M',M}=\Theta\circ\tilde f\circ P:\, M'\otimes M\lra M\otimes M'$$
is an isomorphism of $U_{r,s}(\mathfrak g)$-modules, where
$P:\,M'\otimes M\lra M\otimes M'$ is the flip map such that
$P(m'\otimes m)=m\otimes m'$ for any $m\in M,\,m'\in M'$.
\end{theorem}
\begin{proof} \ Obviously, $R_{M',M}$ is invertible. It remains to show
that $R_{M',M}$ is a $U_{r,s}(\mathfrak g)$-module homomorphism,
that is, to check that
$$
\Delta(a)R_{M',M}(m'\otimes m)=R_{M',M}\Delta(a)(m'\otimes m)
\leqno(9)
$$
holds for all $a\in U_{r,s}(\mathfrak g)$, $m\in M_\lambda$ and
$m'\in M_\mu'$. It suffices to verify (9) for the generators
$e_n,\,f_n,\,\om_n,\,\om_n'$, because the subalgebra generated by
the first $4(n-1)$ generators $e_i,\,f_i,\,\om_i,\,\om_i'$ $(1\le
i<n)$ is isomorphic to $U_{r,s}(\mathfrak{sl}_n)$, and this can be
reduced to the proof of the type $A$ case (see [BW1, Theorem 4.11]).
We will present the computation just for $a=f_n$. Using Lemma 3.3
(iii), we get
\begin{equation*}
\begin{split}
\text{LHS of (9)}&=f(\lambda,\mu)\Delta(f_n)\Theta(m\otimes
m')\\
&=f(\lambda,\mu)(1\otimes
f_n)\bigl(\sum\Theta_\zeta\bigr)(m\otimes m')\\
&\qquad+f(\lambda,\mu)(f_n\otimes
\om_n')\bigl(\sum\Theta_{\zeta-\al_n}\bigr)(m\otimes m')\\
&=f(\lambda,\mu)\bigl(\sum\Theta_\zeta\bigr)(1\otimes
f_n)(m\otimes m')\\
&\qquad+f(\lambda,\mu)\bigl(\sum\Theta_{\zeta-\al_n}\bigr)(f_n\otimes
\om_n)(m\otimes m')\\
&=f(\lambda,\mu)\lg\om_n',\om_\lambda\rg\bigl(\sum\Theta_\zeta\bigr)(\om_n'm\otimes f_nm')\\
&\qquad+f(\lambda,\mu)\lg
\om_\mu',\om_n\rg\bigl(\sum\Theta_{\zeta-\al_n}\bigr)(f_nm\otimes
m').
\end{split}
\end{equation*}
On the other hand, we have
\begin{equation*}
\begin{split}
\text{RHS of (9)}&=R_{M',M}(m'\otimes f_nm+f_nm'\otimes
\om_n'm)\\
&=(\Theta\circ\tilde f)(f_nm\otimes m'+\om_n'm\otimes f_nm')\\
&=f(\lambda-\al_n,\mu)\Theta(f_nm\otimes
m')+f(\lambda,\mu-\al_n)\Theta(\om_n'm\otimes f_nm')\\
&=f(\lambda-\al_n,\mu)\bigl(\sum\Theta_\zeta\bigr)(f_n\otimes
1)(m\otimes m')\\
&\qquad+f(\lambda,\mu-\al_n)\bigl(\sum\Theta_\zeta\bigr)(\om_n'\otimes
f_n)(m\otimes m').
\end{split}
\end{equation*}
Thus (3) indicates that (9) holds.
\end{proof}

\begin{remark} \ Similar to the treatment in [BW1,
Section 5] for the type $A$ case, we can prove the maps $R_{M',M}$
satisfy the quantum Yang-Baxter equation for our cases. That is,
given three $U_{r,s}(\mathfrak g)$-modules $M,\,M',\,M''$ in
$\mathcal O$, we have $R_{12}\circ R_{23}\circ R_{12}=R_{23}\circ
R_{12}\circ R_{23}$ as maps from $M\otimes M'\otimes M''$ to
$M''\otimes M'\otimes M$ (see [BW1, Theorem 5.4]). On the other
hand, we also can prove the hexagon identities (see [BW1, Theorem
5.7]) for the maps $R_{M',M}$ by the same approach. Consequently,
$\mathcal O$ is a braided monoidal category with braiding
$R=R_{M',M}$ for each pair of modules $M',\,M$ in $\mathcal O$.
\end{remark}

\subsection*{Quantum Casimir operators and complete reducibility.}
\smallskip

The $U_{r,s}(\mathfrak g)$-module isomorphisms $R_{M',M}$
constructed in Theorem 3.4, which are called the $R$-matrices, are
mainly determined by $\Theta$. For the expression (7) of
$\Theta_\zeta$, we set
\begin{gather*}
\Omega_\zeta=\sum_kS(v_k^{\zeta})u_k^\zeta, \qquad
\Omega_\zeta'=\theta(\Omega_\zeta),\tag{10}\\
\Omega=\sum_{\zeta\in Q^+}\Omega_\zeta,\qquad
\Omega'=\sum_{\zeta\in Q^+}\Omega_\zeta',\tag{11}
\end{gather*}
where $\theta$ is the $\mathbb Q$-algebra isomorphism of
$U_{r,s}(\mathfrak g)$ into its associated quantum group
$U_{s^{-1},r^{-1}}(\mathfrak g)$ (for definition, see [BGH])
introduced in the formula (11) in Section 2. Obviously,
$\Theta_\zeta$, $\Omega_\zeta$, $\Omega$ and $\Omega'$ are
independent of the choice of bases $\{u_k^\zeta\}$ and
$\{v_k^\zeta\}$. $\Omega$ preserves the weight spaces of any
$M\in\mathcal O$.

\begin{definition} \ {\it The element $\Omega$ is called a quantum
Casimir element for the two-parameter quantum group
$U_{r,s}(\mathfrak g)$.}
\end{definition}

\begin{prop} \ Let $\psi$ and $\varphi$ be the algebra
automorphisms of $\,U_{r,s}(\mathfrak g)$ such that
$\psi(\om_i)=\om_i$, $\psi(\om_i')=\om_i'$,
$\psi(e_i)=\om_i'\om_i^{-1}e_i$, $\psi(f_i)=f_i{\om_i'}^{-1}\om_i$
and $\varphi(\om_i)=\om_i,\,\varphi(\om_i')=\om_i'$,
$\varphi(e_i)=e_i\om_i^{-1}\om_i'$,
$\varphi(f_i)=\om_i{\om_i'}^{-1}f_i$. Then
$$
\psi(a)\,\Omega=\Omega\,a, \qquad \varphi(a)\,\Omega'=\Omega'\,a,
\quad\text{\it for }\; a\in U_{r,s}(\mathfrak g).\leqno(12)
$$
\end{prop}
\begin{proof} \ Since $\psi$ is an algebra automorphism, it is
enough to prove the first assertion for the generators
$a=\om_i,\,\om_i',\,e_i,\,f_i$. For $a=\om_i$ or $\om_i'$, it is
obviously true. Applying the mapping $\mathfrak m\circ(S\otimes
1)$ to both sides of Lemma 3.3 (ii) \& (iii) (where $\mathfrak m$
is the product of $U_{r,s}(\mathfrak g)$ and $S$ is its antipode)
and summing over $\zeta\in Q^+$ we obtain
$\Omega\,e_i=\om_i'\om_i^{-1}e_i\Omega$ and $\Omega\,f_i
=f_i{\om_i'}^{-1}\om_i\Omega$. This means that
$\Omega\,e_i=\psi(e_i)\,\Omega$ and
$\Omega\,f_i=\psi(f_i)\,\Omega$. Applying the automorphism
$\theta$ we get the assertion for $\Omega'$.
\end{proof}

\begin{coro} \ {\it For $M\in\text{\rm Ob}(\mathcal O)$, assume that
$m\in M_\lambda$. Then}

$(\text{\rm i})$ \ \
$\Omega\,e_i.m=(r_is_i^{-1})^{-(\lambda+\al_i,\al_i^\vee)}e_i\,\Omega.m,$

$(\text{\rm ii})$ \
$\Omega\,f_i.m=(r_is_i^{-1})^{(\lambda,\al_i^\vee)}f_i\,\Omega.m.$
\end{coro}
\begin{proof} \ By Proposition 3.7, for $m\in M_\lambda$, we have
\begin{gather*}
\psi(e_i)\,\Omega.m=\lg
\om_{\lambda+\al_i}',\om_i\rg^{-1}\lg\om_i',\om_{\lambda+\al_i}\rg^{-1}e_i\,\Omega.m,\\
\psi(f_i)\,\Omega.m=\lg\om_\lambda',\om_i\rg\,\lg\om_i',\om_\lambda\rg\,f_i\,\Omega.m.
\end{gather*}
Using formulas (B), (C), \& (D) in Lemma 2.1, we can conclude the
required result.
\end{proof}

\begin{remark} \ According to Section 1, we have made a convention:
we have $r_i=r^{(\al_i,\al_i)},\,s_i=s^{(\al_i,\al_i)}$ only in the
type $B$ case,so
$(r_is_i^{-1})^{(\lambda,\al_i^\vee)}=(rs^{-1})^{2(\lambda,\al_i)}$
for any $i$. However, for any other case, we always have
$(r_is_i^{-1})^{(\lambda,\al_i^\vee)}=(rs^{-1})^{(\lambda,\al_i)}$
for any $i$ since
$r_i=r^{\frac{(\al_i,\al_i)}{2}},\,s_i=s^{\frac{(\al_i,\al_i)}{2}}$.
Based on this observation, we make the following definition.
\end{remark}

\begin{definition} \ {\it For $M\in\text{\rm Ob}(\mathcal O)$,
define a linear operator $\omega:\, M\lra M$ by setting
$$
\omega.v_\mu=(rs^{-1})^{\frac{\Delta_{X,B}}{2}(\mu+\rho,\mu+\rho)}v_\mu,
\quad \text{\it for } \ v_\mu\in M_\mu,\leqno(13)
$$
where $\rho$ is the half-sum of all positive roots of $\mathfrak
g$, and $\Delta_{X,B}=2$ if $X=B$, otherwise, $\Delta_{X,B}$ will
take value $1$.}
\end{definition}

\begin{prop} \ Assume that the Verma module
$M(\lambda)\in\text{\rm Ob}(\mathcal O)$, then the operator
$\Omega\omega$ is a multiple of the identity operator, that is,
$$
\Omega\omega=(rs^{-1})^{\frac{\Delta_{X,B}}2(\lambda+\rho,\lambda+\rho)}I.\leqno(14)
$$
\end{prop}
\begin{proof} \ Let $v_\lambda$ be a highest weight vector of the
Verma module $M(\lambda)$. Then $M(\lambda)=U_{r,s}(\mathfrak
n^-)v_\lambda=\sum_{\be\in Q^+}U_{r,s}^{-\be}(\mathfrak
n^-)v_\lambda$. For $f_\be\in U_{r,s}^{-\be}(\mathfrak n^-)$,
denote $v_{\lambda-\be}:=f_\be.v_\lambda$, which is a weight
vector of weight $\lambda-\be$. We claim that
$$
\Omega\omega.f_i.v_{\lambda-\be}=f_i.\Omega\omega.v_{\lambda-\be},\leqno(15)
$$
for any $\be\in Q^+$ and any $i$. Indeed, noting that
$$
\frac1{2}\bigl[\,(\si-\al_i+\rho,\si-\al_i+\rho)-(\si+\rho,\si+\rho)\,\bigr]+(\si,\al_i)=\frac1{2}\bigl[\,(\al_i,\al_i)
-2(\al_i,\rho)\,\bigr]=0,
$$
and setting $\lambda-\be=\si$, we have
\begin{equation*}
\begin{split}
 \Omega\omega.f_i.v_{\lambda-\be}&=(\Omega
f_i)(rs^{-1})^{\Delta_{X,B}c}\,\omega.v_{\lambda-\be}\\
&=(f_i{\om_i'}^{-1}\om_i\Omega)\,(rs^{-1})^{\Delta_{X,B}c}\,\omega.v_{\lambda-\be}\\
&=f_i(rs^{-1})^{\Delta_{X,B}(\lambda-\be,\al_i)}(rs^{-1})^{\Delta_{X,B}c}\,\Omega\omega.v_{\lambda-\be}\\
&=f_i\,\Omega\omega.v_{\lambda-\be},
\end{split}
\end{equation*}
where
$c=\frac1{2}\bigl[\,(\lambda-\be-\al_i+\rho,\lambda-\be-\al_i+\rho)-(\lambda-\be+\rho,\lambda-\be+\rho)\,\bigr]$.
(15) yields
\begin{equation*}
\begin{split}
\Omega\omega.f_\be.v_\lambda&=f_\be.\Omega\omega.v_\lambda\\
&=(rs^{-1})^{\frac{\Delta_{X,B}}2(\lambda+\rho,\lambda+\rho)}f_\be.\Omega.v_\lambda\\
&=(rs^{-1})^{\frac{\Delta_{X,B}}2(\lambda+\rho,\lambda+\rho)}f_\be.\Omega_0.e_\lambda\\
&=(rs^{-1})^{\frac{\Delta_{X,B}}2(\lambda+\rho,\lambda+\rho)}f_\be.v_\lambda,\quad
(\Omega_0=1).
\end{split}
\end{equation*}
So the relation (14) follows.
\end{proof}

\begin{coro} \ $(\text{\rm i})$ \ {\it For the simple
$U_{r,s}(\mathfrak g)$-module $L(\lambda)\in\text{\rm Ob}(\mathcal
O)$, there holds}
$$\Omega\omega=(rs^{-1})^{\frac{\Delta_{X,B}}2(\lambda+\rho,\lambda+\rho)}I.$$
$(\text{\rm ii})$ \ {\it For each finite-dimensional
$M\in\text{\rm Ob}(\mathcal O)$, the eigenvalues of the operator
$(\Omega\omega)|_M$ are integral powers of
$(rs^{-1})^{\frac1{2}}$.}
\end{coro}
\begin{proof} \ (i) is evident. For (ii), as $M\in\text{\rm
Ob}(\mathcal O)$ is finite-dimensional, it has a composition
series whose factors are finite-dimensional simple
$U_{r,s}(\mathfrak g)$-modules in $\mathcal O$, on which
$\Omega\omega$ acts as multiplication by
$(rs^{-1})^{\frac{\Delta_{X,B}}2(\mu+\rho, \mu+\rho)}$ for some
$\mu\in\Lambda^+$, as indicated by (i) and Corollary 2.8. After
taking an appropriate basis of $M$ compatible with a chosen
composition series, the acting matrix of $(\Omega\omega)|_M$ has
the required property.
\end{proof}

From Corollary 3.8 and Definition 3.10, we have a further result as
follows.

\begin{theorem} \ The operator $\Omega\omega:\,M\lra
M$ commutes with the action of $U_{r,s}(\mathfrak g)$ on any
module $M\in\text{\rm Ob}(\mathcal O)$, where $\mathfrak
g=\mathfrak
{sl}_{n+1},\,\mathfrak{so}_{2n+1},\,\mathfrak{so}_{2n}$, or
$\mathfrak {sp}_{2n}$.
\end{theorem}
\begin{proof} \ At first, it needs to show that $\Omega\omega$
commutes with $e_i$, $f_i$ $(1\le i\le n)$. For $m\in M_\mu$, by
Corollary 3.8 and Definition 3.10, we get
\begin{equation*}
\begin{split}
\Omega\omega.(e_i.m)&=(rs^{-1})^{\frac{\Delta_{X,B}}2(\mu+\al_i+\rho,\,\mu+\al_i+\rho)}\Omega
e_i.m\\
&=(rs^{-1})^{\Delta_{X,B}[\,\frac{1}2(\mu+\al_i+\rho,\,\mu+\al_i+\rho)-(\mu+\al_i,\,\al_i)\,]}e_i\Omega.m\\
&=(rs^{-1})^{\frac{\Delta_{X,B}}2(\mu+\rho,\,\mu+\rho)}e_i\Omega.m\\
&=e_i.(\Omega\omega.m).\\
\Omega\omega.(f_i.m)&=(rs^{-1})^{\frac{\Delta_{X,B}}2(\mu-\al_i+\rho,\,\mu-\al_i+\rho)}\Omega
f_i.m\\
&=(rs^{-1})^{\Delta_{X,B}[\,\frac{1}2(\mu-\al_i+\rho,\,\mu-\al_i+\rho)+(\mu,\,\al_i)\,]}f_i\Omega.m\\
&=(rs^{-1})^{\frac{\Delta_{X,B}}2(\mu+\rho,\,\mu+\rho)}f_i\Omega.m\\
&=f_i.(\Omega\omega.m).
\end{split}
\end{equation*}
Obviously, $\Omega\omega$ commutes with the action of $\om_i$,
$\om_i'$ $(1\le i\le n)$, for it preserves the weight spaces of
$M$.
\end{proof}

The following Lemma is due to [BW2, Lemma 3.7] for the case of
$\mathfrak g=\mathfrak {gl}_{n+1}$, or $\mathfrak {sl}_{n+1}$,
which is still valid in our cases.

\begin{lemma} \ Assume that $rs^{-1}$ is not a root of
unity. Let $M$ be a nonzero finite-dimensional quotient of the
Verma $U_{r,s}(\mathfrak g)$-module $M(\lambda)\in\text{\rm
Ob}(\mathcal O)$. Then $M$ is simple. In particular,
$L'(\lambda)=L(\lambda)$ for $\lambda\in \Lambda^+$.
\end{lemma}
\begin{proof} \ Lemma 2.6 means $\lambda\in\Lambda^+$. The proof is
based on the counter-evidence method and Proposition 3.11, which is
the same as that of [BW2, Lemma 3.7], with slight differences: for
the function $g(\lambda)$ used in the proof there we use
$(rs^{-1})^{\frac{\Delta_{X,B}}2(\lambda+\rho,\lambda+\rho)}$
instead, noting the fact from Lie algebra theory (see [D], or [K])
that for any weight $\mu\le \lambda$ where $\lambda\in\Lambda^+$,
$(\lambda+\rho, \lambda+\rho)=(\mu+\rho, \mu+\rho)$ if and only if
$\mu=\lambda$.
\end{proof}

Based on the above results, using a similar argument due to Kac
[K] in the proof of complete reducibility of category $\mathcal O$
for affine Kac-Moody Lie algebras (or comparing with the proof of
[BW2, Theorem 3.8] in the spirit of Lusztig [L1]), we have

\begin{theorem} \ Assume that $rs^{-1}$ is a non-root of
unity. For $\mathfrak g=\mathfrak {sl}_{n+1}$, $\mathfrak
{so}_{2n+1}$, $\mathfrak{so}_{2n}$ or $\mathfrak {sp}_{2n}$, let $M$
be a nonzero finite-dimensional $U_{r,s}(\mathfrak g)$-module on
which $U^0$ acts semisimply. Then $M$ is completely reducible.
\hfill\qed
\end{theorem}

\bigskip
\section*{Acknowledgments}
\medskip

Part of this work was done when N. Hu visited the Department of
Mathematics and Statistics of York University in Canada from Sept.
2003 to Sept. 2004, visited the DMA, l'Ecole Normale Sup\'erieure de
Paris in France as an invited professor from Oct. 20 to Nov. 20,
2004, and visited the Fachbereich Mathematik der Universit\"at
Hamburg in Germany as a DFG-visiting professor from Nov. 21, 2004 to
Feb. 28, 2005. He would like to express his deep thanks to M. Rosso
and H. Strade for their invitation and extreme hospitality, as well
as the supports from York University, l'Ecole Normale Sup\'erieure
de Paris and die Deutche Forschungsgemeinschaft (DFG).

\bigskip
\bibliographystyle{amsalpha}

\end{document}